\documentclass[11pt]{article}%
\usepackage{amssymb}
\usepackage[all]{xy}
\usepackage[latin1]{inputenc}
\usepackage{amsfonts}
\usepackage{amsmath}
\usepackage{amssymb}
\usepackage{url}
\usepackage{hyperref}%
\usepackage{graphicx}
\providecommand{\U}[1]{\protect\rule{.1in}{.1in}}
\newtheorem{theo}{Theorem}[section]
\newtheorem{prop}[theo]{Proposition}
\newtheorem{lem}[theo]{Lemma}
\newtheorem{cor}[theo]{Corollary}

\newtheorem{defi}[theo]{Definition}

\numberwithin{equation}{section}

\usepackage{dsfont}
\newcommand{\indicator}[1]{\mathds{1}_{ {#1}}}
\def\un{\indicator{}}

\def\d{\delta}

\def\e{\varepsilon}

\newcommand{\EE}{\mathbb{E}}

\newcommand{\LL}{\mathbb{L}}

\newcommand{\PP}{\mathbb{P}}

\newcommand{\RR}{\mathbb{R}}

\newcommand{\Aa}{ {\cal A }}
\newcommand{\Ba}{ {\cal B }}

\newcommand{\Mb}{ {\bf M }}

\newcommand{\Ca}{ {\cal C }}
\newcommand{\Da}{ {\cal D }}
\newcommand{\La}{ {\cal L }}
\newcommand{\Na}{ {\cal N }}

\newcommand{\Ea}{ {\cal E }}
\newcommand{\Sa}{ {\cal S }}
\newcommand{\Ra}{ {\cal R }}

\newcommand{\Ua}{ {\cal U }}
\newcommand{\Fa}{ {\cal F }}
\newcommand{\Ga}{ {\cal G }}

\newcommand{\Ia}{ {\cal I }}
\newcommand{\Xa}{ {\cal X }}
\newcommand{\Ma}{ {\cal M }}
\newcommand{\Ya}{ {\cal Y}}
\newcommand{\Ha}{ {\cal H }}
\newcommand{\Ja}{ {\cal J }}
\newcommand{\Pa}{ {\cal P }}
\newcommand{\Za}{ {\cal Z }}

\newcommand{\point}{\mbox{\LARGE .}}
\newcommand{\proof}{\noindent\mbox{\bf Proof:}\\}
\newcommand{\cqfd}{\hfill\blbx \\}
\def\blbx{\hbox{\vrule height 5pt width 5pt depth 0pt}\medskip}

\def \PP{\mathbb{P}}
\def \RR{\mathbb{R}}
\def \EE{\mathbb{E}}

\def \LL{\mathbb{L}}

\def \e{\epsilon}
\begin{document}

\title{Moderate Deviations for Mean Field Particle Models}
\author{Pierre Del Moral\thanks{Centre INRIA Bordeaux Sud-Ouest Institut de Math\'ematiques,
Universit\'e Bordeaux I,33405 Talence cedex, France.\
\texttt{pierre.del-moral@inria.fr.}},
 Shulan Hu\thanks{Corresponding author: Department of Statistics, Zhongnan University of Economics
and Law, Wuhan {\rm 430073}, China. \
 \texttt{hu\_shulan@yahoo.com}.
}, Liming Wu\thanks{Laboratoire de Math\'ematiques Appliqu\'ees,
CNRS-UMR 6620, Universit\'e Blaise Pascal, 63177 Aubiere, France.\
\texttt{li-ming.wu@math.univ-bpclermont.fr}.\ \ Institute of Applied
Mathematics, Chinese Academy of Sciences, 100190 Beijing, China.\
\texttt{wuliming@amt.ac.cn}.}}

 \maketitle

\begin{abstract}
This article is concerned with  moderate deviation principles of a
general class of mean field type interacting particle models. We
discuss functional moderate deviations of the occupation measures
for both the strong $\tau$-topology on the space of finite and
bounded measures as well as for the corresponding stochastic
processes on some class of functions equipped with the uniform
topology. Our approach is based on an original semigroup analysis
combined with stochastic perturbation techniques and projective
limit large deviation methods.



\end{abstract}


\section{Introduction}
\subsection{Mean Field Particle Models}
Let $(E_n)_{n\geq 0}$ be a sequence of measurable spaces equipped with some $\sigma$-fields $(\Ea_n)_{n\geq 0}$,
 and we let $\Pa(E_n)$ be the set of all probability measures over the set $E_n$,
with $n\geq 0$. We consider a collection
of transformations $\Phi_n~:~\Pa(E_{n-1})\rightarrow\Pa(E_n)$ and we denote by $(\eta_n)_{n\geq 0}$
 a sequence of probability measures  on $E_n$ that satisfies a nonlinear
equation of the following form
\begin{equation}\label{phi}
\eta_{n+1}=\Phi_{n+1}\left(\eta_{n}\right).
\end{equation}

 The mean field  particle interpretations of  these measure valued models
relies on the fact that the one step mappings can be rewritten in the following form
\begin{equation}\label{phidef}
\Phi_n\left(\eta_{n-1}\right)=\eta_{n-1}K_{n,\eta_{n-1}}
\end{equation}
for some collection of Markov kernels $K_{n,\mu_{n-1}}$ indexed by
the time parameter $n$ and the set of probability measures $\mu_n$
on the space $E_{n-1}$.  These models provide a natural
interpretation of the distribution laws $\eta_n$ as the laws of a
non linear Markov chain whose elementary transitions depend on the
current distribution.  In the further development of the article, we
always assume that  the mappings
$$
\left(x_n^i\right)_{1\leq i\leq N}\in E^N_n\mapsto K_{n+1,\frac{1}{N}\sum_{j=1}^N\delta_{x^j_n}}\left(x^i_n,A_{n+1}\right)
$$
are $\Ea^{\otimes N}_n$-measurable, for any $n\geq 0$,  $1\leq i\leq N$, and any measurable subset $A_{n+1}\subset E_{n+1}$.
In this situation,
the mean field particle interpretation of this nonlinear measure valued model is an $E^N_n$-valued Markov chain \index{$\xi^{(N)}_n$}
$
\xi^{(N)}_n=\left(\xi^{(N,i)}_n\right)_{1\leq i\leq N}
$,
with elementary transitions
defined as
\begin{equation}\label{meanfield}
\PP\left(\xi^{(N)}_{n+1}\in dx
  ~\left|~\Aa^{(N)}_n\right.\right)
=\prod_{i=1}^N~K_{n+1,\eta^N_n}( \xi^{(N,i)}_n,dx^i) \quad\mbox{\rm
with}\quad \eta^N_n:=\frac{1}{N}\sum_{j=1}^N\delta_{\xi_n^{(N,j)}}.
\end{equation}
In the above displayed formula,
$\Aa^{(N)}_n:=\sigma\left(\xi^{(N)}_p,~0\leq p\leq n\right)$ stands
for the sigma-field generated by the random variables
$(\xi^{(N)}_p)_{0\leq p\leq n}$, and $dx=dx^1\times\ldots\times
dx^N$ stands for an infinitesimal neighborhood of a point
$x=(x^1,\ldots,x^N)\in E_n^N$. The initial system $\xi^{(N)}_0$
consists of $N$ independent and identically distributed random
variables with common law $\eta_0$. To simplify the presentation,
when there is no possible confusion we suppress the parameter $N$,
so that we write $\xi_n$ and $\xi^i_n$ instead of $\xi^{(N)}_n$ and
$\xi^{(N,i)}_n$. For a thorough description of these discrete
generation and non linear McKean type models, we refer the reader
to~\cite{fk}.

 A
typical example we have in mind is the Feynman-Kac model associated
with $(0,1]$-valued potential functions $G_n$ and Markov transitions
$M_{n+1}$ from $E_n$ into $E_{n+1}$ given by
\begin{equation}\label{phiFK}
\Phi_{n+1}\left(\eta_{n}\right)(dy)=\left(\Psi_{G_{n}}\left(\eta_n\right)M_{n+1}\right)(dy):=
\int~\Psi_{G_{n}}\left(\eta_n\right)(dx)~M_{n+1}(x,dy)
\end{equation}
where $\displaystyle \Psi_G(\eta)(dx):=\frac{G(x)}{\eta(G)}\eta(dx)$
(with $\eta(G):=\int Gd\eta(x)$).

In this situation, the flow of measures $\eta_n$ is given for any
bounded measurable function $f$ on $E_n$ by the following functional
formula
$$
\eta_n(f_n)=\int_{E_n} f_n(x)~\eta_n(dx)\propto
\EE\left(f_n(X_n)~\prod_{0\leq p<n}G_p(X_p)\right)
$$
where $X_n$ stands for a Markov chain with initial distribution
$\eta_0$ and Markov transitions $M_n$.

 Recall that $\Psi_{G_n}(\eta_n)$ can be
expressed as a non-linear Markov transport equation
\begin{equation}\label{nonlin-transport}
\Psi_{G_n}(\eta_n)=\eta_n S_{\eta_n,G_n}
\end{equation}
with the Markov transitions
$$
S_{\eta_n,G_n}(x,dy)=G_n(x)~\delta_x(dy)+\left(1-G_n(x)\right)~\Psi_{G_n}(\eta_n)(dy)
$$
we find that
$$
K_{n+1,\eta_{n}}=S_{\eta_{n},G_{n}}M_{n+1}
$$

These measure valued equations arise in a variety of applications areas,
including in physics, biology and in advanced stochastic engineering sciences. For instance,
 in signal processing, the
conditional distributions of  the paths  of Markov signal  given a
series of noisy observations satisfy a two-step prediction-updating
equation of the form (\ref{phi}).  In this context, the state space
$E_n$ depends on the time parameter and it consists of all signal
path sequences all length $n$. In this situation, it is worth
mentioning that the corresponding mean field particle model in path
space represents the evolution of a genealogical tree model
associated with a genetic type algorithm.

In the context of sequential bayesian inference, the distributions
$\eta_n$ could also be the posterior distributions of an unknown
parameter given  the data collected up to time $n$. These equations
also arise in physics and in molecular chemistry. In this situation,
$\eta_n$ is often interpreted as the law of a particle  evolving in
an absorbing medium related to some potential energy function. These
non linear models are also used in advanced stochastic engineering
sciences, and more particularly in stochastic optimization as well
as in rare event simulation. In these situations, $\eta_n$ is often
given by a Boltzmann-Gibbs measure associated with some decreasing
temperature parameter or some decreasing sequence of critical rare
event levels. In the late case, the state spaces $E_n$ represent the
set of excursions of the reference Markov chain between two level
sets.

During the last two decades, the mean field particle interpretations
of these discrete generation measure valued equations are
increasingly identified as a powerful stochastic simulation
algorithm. They have led to spectacular results in signal processing
with the corresponding  particle filter technology, in stochastic
engineering with interacting type Metropolis and Gibbs sampler
methods, as well as in quantum chemistry with quantum and diffusion
Monte Carlo algorithms leading to precise estimates of the top
eigenvalues and the ground states of Schroedinger operators.  It is
clearly out of the scope of this article to review these models. For
a more detailed discussion on these application areas, we again
refer the reader to~\cite{fk,ddj,arnaud,dmrio}, and the references
therein.

The mathematical and numerical analysis of these mean field particle
models (\ref{meanfield}) is one of the most active research subject in pure and applied probability,
as well as in advanced
stochastic engineering and computational physics.
In recent years, a variety of mathematical results have been discussed in the literature,
including propagation of chaos type properties, $\LL_p$-mean error bounds, as well as fluctuations theorems, large deviation principles
and non asymptotic concentration inequalities.
The moderate deviation
properties can be thought as an intermediate asymptotic estimation between the central limit theorem and the large deviations principles.
 Most of the existing literature on
moderate deviation principles is concerned with independent and identically distributed
random sequences or  Markov chain processes;
 see for instance the series of works by M.A. Arcones~\cite{arcones, arcones2},
  A. de Acosta~\cite{deacosta},  A. de Acosta and X. Chen~\cite{deacostachen},
H. Djellout and A. Guillin~\cite{guillin1},
F.Q. Gao~\cite{gao1,gao2}, M. Ledoux~\cite{ledoux},
L.M. Wu~\cite{liming,liming2}.

Surprisingly very few articles discuss moderate deviations
for mean field interacting particle models.
The first pioneering article discussing moderate deviations for interacting
processes seems to be the work by R. Douc, A. Guillin and J. Najim~\cite{douc}. In this  article,
the authors prove a moderate deviation for the empirical mean  value of a particle filter model associated
with some bounded and unbounded fixed sequence of test functions. In our framework, we also mention that the particle
filter stochastic model discussed in~\cite{douc} is associated with a class of McKean type
transitions of the form $K_{n,\eta}(x,dy)=\Phi_n(\eta)(dy)$. The main simplification due to this
choice of transition comes from the fact that the corresponding mean field particle model reduces
to a collection of conditionally independent and identically distributed random variables.

The rather weak regularity properties used in this work  follow a
recent article of the first author with E. Rio~\cite{dmrio}.
 The mathematical framework developed in this recent work
applies to a general class of
mean field particle models, including  Feynman-Kac integration models,
 McKean Vlasov diffusion type models, as well as  McKean collision type models of gases.
We generalized the classical Hoeffding, Bernstein and Bennett
inequalities for independent random sequences to interacting
particle systems but we left open the question of moderate deviation
principles. In the present article, we complete this study with
functional moderate deviations of mean field particle models
 for both the
$\tau$-topology on the space of signed and bounded measures and for the empirical random field processes
associated with some collection of functions.  Our analysis is based on
 an original semigroup analysis  combined with stochastic
perturbation techniques and projective limit deviation methods.

 \subsection{Outline of the paper.} This paper is organized as
 follows. In the next section we present the main results : the moderate deviation principles (MDP in short) in three types : (1) in finite
 dimension ;
 (2) in infinite dimension but for the $\tau$-topology; (3) for
 empirical process indexed by a class of functions; and we describe some main
 lines leading to them. We prove the MDP in finite dimension in
 section 3. We prove in Section 4 the MDP in
 the $\tau-$topology by the method of projective limit. We establish in Section 5 the MDP for
 empirical processes by the method of metric entropy. Some
 complicated and technical results needed in the proofs of the main
 theorems are provided in the three Appendices : Section 6, 7
 and 8.

\subsection{Some notations}
We end this introduction with some more or less traditional
notations used in the present article. We denote respectively by
$\mathcal{M}(E)$, $\mathcal{M}_{0}(E)$, and $\mathcal{B}(E)$, the
set of all finite signed measures on some measurable space
$(E,\mathcal{E})$, the convex subset of finite signed measures $\nu$
with $\nu(E)=0$, and the Banach space of all bounded and measurable
functions $f$ equipped  with the uniform norm $\Vert f\Vert$. We
also denote by
 $\mbox{Osc}_{1}(E)$, the
convex set of $\mathcal{E}$-measurable functions $f$ with
oscillations $\mbox{osc}(f):=\sup_{x\not=y}|f(x)-f(y)|\leq1$. We let
$\mu(f)=\int~\mu(dx)~f(x)$, be the Lebesgue integral of a function
$f\in\mathcal{B}(E)$, with respect to a measure
$\mu\in\mathcal{M}(E)$. We recall that a bounded integral operator
$M$ from a measurable space $(E,\mathcal{E})$ into an auxiliary
measurable space $(F,\mathcal{F})$ is an operator $f\mapsto M(f)$
from $\mathcal{B}(F)$ into $\mathcal{B}(E)$ such that the functions
$x\mapsto M(f)(x):=\int_{F}M(x,dy)f(y)$ are $\mathcal{E}$-measurable
and bounded, for any $f\in\mathcal{B}(F)$. A Markov kernel is a
positive and bounded integral operator $M$ with $M(1)=1$. Given a
pair of bounded integral operators $(M_1,M_2)$, we let $(M_1M_2)$
the composition operator defined by $(M_1M_2)(f)=M_1(M_2(f))$. For
time homogenous state spaces, we denote by $M^m=M^{m-1}M=MM^{m-1}$
the $m$-th composition of a given bounded integral operator $M$,
with $m\geq 1$.

A
bounded integral operator $M$ from a measurable space $(E,\mathcal{E})$ into
an auxiliary measurable space $(F,\mathcal{F})$ also generates a dual operator
$\mu\mapsto\mu M$ from $\mathcal{M}(E)$ into $\mathcal{M}(F)$ defined by $(\mu
M)(f):=\mu(M(f))$. We let  $b(m)$ be the collection of constants  given below
$$
b(2m)^{2m}:=\frac{(2m)!}{m! 2^m},\qquad\mbox{and}\qquad
b(2m+1)^{2m+1}:=\frac{(2m+1)!}{(m+1)!\sqrt{m+1/2}}~2^{-(m+1/2)}.
$$
 When the bounded integral
operator $M$ has a constant mass, that is $M(1)\left( x\right)
=M(1)\left( y\right) $ for any $(x,y)\in E^{2}$, the operator
$\mu\mapsto\mu M$ maps $\mathcal{M}_{0}(E)$ into
$\mathcal{M}_{0}(F)$. In this situation, we let $\beta(M)$ be the
Dobrushin coefficient of a bounded integral operator $M$ defined by
the following formula \begin{equation}\label{beta} \beta(M):=\sup{\
\{\mbox{\rm osc}(M(f))\;;\;\;f\in\mbox{\rm
Osc}_{1}(F)\}}.\end{equation}

Finally, we  let $\Phi_{p,n}$, $0\leq p\leq n$, be the semigroup
associated with the measure valued equation defined in (\ref{phi}).
 $$
 \Phi_{p,n}=\Phi_n\circ\Phi_{n-1}\circ\ldots\circ\Phi_{p+1}.
 $$
 For $p=n$, we use the convention $\Phi_{n,n}=Id$, the identity  operator.


\section{Description of the main results and a first order fluctuation analysis}
\subsection{Regularity conditions}\label{regsec}
We let $\Upsilon(E_1,E_2)$ be the set of mappings $\Phi~:~\mu\in \Pa(E_1)\mapsto \Phi(\mu)\in \Pa(E_2)$ satisfying the first order decomposition
\begin{equation}\label{firstodef}
 \Phi(\mu)-\Phi(\eta)=(\mu-\eta) D_{\eta}\Phi+\Ra^{\Phi}(\mu,\eta)
 \end{equation}
where \begin{description} \item{(i)} the first order operators
$(\Da_{\eta}\Phi)_{\eta\in\Pa(E_1)}$  is some collection of bounded
integral operators  from $E_1$ into $E_2$  such that $\forall
\eta\in \Pa(E_1),~ \forall x\in E_1, (D_{\eta}\Phi)(1)(x)=0$ and
\begin{equation}\label{firsto}\beta
\left(\Da\Phi\right):=\sup_{\eta\in \Pa(E_1)}\beta
\left(D_{\eta}\Phi\right)<\infty;
\end{equation}
\item{(ii)} the collection of second order remainder signed measures
$(\Ra^{\Phi}(\mu,\eta))_{(\mu,\eta)\in\Pa(E^2_1)}$ on $E_2$ are such
that
\begin{equation}\label{condxi222}\left|\Ra^{\Phi}(\mu,\eta)(f)\right|\leq \int~\left|(\mu-\eta)^{\otimes
2}(g)\right|~R^{\Phi}_{\eta}(f,dg)
\end{equation}
for some collection of integral operators $R^{\Phi}_{\eta}$ from $\Ba(E_2)$ into the set $\mbox{\rm Osc}_{1}(E_1)^2$ such that
\begin{equation}\label{lipr}
\sup_{\eta\in\Pa(E_1)} \int~\mbox{\rm osc}(g_1)~\mbox{\rm osc}(g_2)~R^{\Phi}_{\eta}(f,d(g_1\otimes g_2))\leq
 ~\mbox{\rm osc}(f)~\delta\left(R^{\Phi}\right)\quad\mbox{\rm with}\quad
 \delta\left(R^{\Phi}\right)<\infty.
\end{equation}
\end{description}

We say that a collection of
 Markov transitions $K_{\eta}$ from a measurable space $(E_1,\Ea_1)$ into another $(E_2,\Ea_2)$
 satisfy condition ($K$) as soon as the
 following Lipschitz type inequality is met  for every $f\in\mbox{\rm Osc}_1(E_2)$:
 \begin{equation}\label{lipsK1}
\hskip-3.6cm(K)\hskip3.3cm \Vert\left[
K_{\mu}-K_{\eta}\right](f)\Vert\leq\int~
|(\mu-\eta)(h)|~T^K_{\eta}(f,dh).
 \end{equation}
 In the above display, $T^K_{\eta}$ stands for some collection of bounded integral operators   from $\Ba(E_2)$ into $\Ba(E_1)$ such that
 \begin{equation}\label{lipsK2}
\sup_{\eta\in \Pa(E_1)} \int~  \mbox{\rm osc}(h)~T^K_{\eta}(f,dh)\leq ~\mbox{\rm osc}(f)~\delta\left(T^{K}\right)
 \end{equation}
for some finite constant $\delta\left(T^{K}\right)<\infty$.  In the
special case where $K_{\eta}(x,dy)=\Phi(\eta)(dy)$, for some mapping
$\Phi~:~\eta\in\Pa(E_1)\mapsto \Phi(\eta)\in \Pa(E_2)$, condition
(K) is a simple Lipschitz type condition on the mapping $\Phi$. In
this situation, we denote by $(\Phi)$ the corresponding condition;
and whenever it is met, we says that the mapping $\Phi$ satisfy
condition $(\Phi)$.

Throughout this paper we assume
\begin{description}

\item{\bf (H1)}
{\em The given collection of  McKean transitions $K_{n,\eta}$
satisfyies the Lipschitz type condition stated in (\ref{lipsK1}) and
(\ref{lipsK2}). We also assume that the one step mappings
$$
\Phi_n ~:~\mu\in\Pa(E_{n-1})\longrightarrow \Phi_n(\mu):=\mu  K_{n,\mu}\in \Pa(E_n)
$$
governing the equation (\ref{phi}) are chosen so that $\Phi_n\in
\Upsilon(E_{n-1},E_n)$, for any $n\geq 1$.} \end{description}

 Several examples of non linear semigroups satisfying these
weak regularity can be found in~\cite{dmrio}, including gaussian type mean field models, and McKean
velocity models of gases. We illustrate our assumptions in the context of Feynman-Kac type models.
In this situation, we have the easily checked formulae
\begin{eqnarray*}
\left[\Phi_{n+1}\left(\mu\right)-\Phi_{n+1}\left(\eta\right)\right](f)&=&\frac{1}{\mu(G_{n,\eta})}~\left(\mu-\eta\right)\left[
G_{n,\eta}~M_{n+1,\eta}(f)\right]\\
&&\\
&=&\left(\mu-\eta\right)\left[
G_{n,\eta}~M_{n+1,\eta}(f)\right]\\
&&\qquad +\frac{1}{\mu(G_{n,\eta})}~\left[\eta-\mu\right](G_{n,\eta})~\left(\mu-\eta\right)\left[
G_{n,\eta}~M_{n+1,\eta}(f)\right]
\end{eqnarray*}
with the functions
$$
G_{n,\eta}=G_n/\eta(G_{n})\quad\mbox{\rm and}\quad M_{n+1,\eta}(f):=M_{n+1}(f)-\Phi_{n+1}\left(\eta\right)(f)
$$
Assuming that $g_n=\sup_{x,y}{G_n(x)/G_n(y)}<\infty$, we find the Lipschitz estimates
\begin{equation}\label{ref-lipsch}
\left|\left[\Phi_{n+1}\left(\eta\right)-\Phi_{n+1}\left(\eta\right)\right](f)\right|\leq g_n~\left|\left(\mu-\eta\right)D_{\eta}\Phi_{n+1}(f)
\right|
\end{equation}
as well as the first order estimation
$$
\left|\left[\left[\Phi_{n+1}\left(\eta\right)-\Phi_{n+1}\left(\eta\right)\right]-\left(\mu-\eta\right)D_{\eta}\Phi_{n+1}\right](f)\right|
\leq g_n~\left|\left[\eta-\mu\right](G_{n,\eta})\right|~\left|\left(\mu-\eta\right)\left[D_{\eta}\Phi_{n+1}(f)\right]\right|
$$
with the first order functional
$$
D_{\eta}\Phi_{n+1}(f)=
G_{n,\eta}~M_{n+1,\eta}(f)
$$

We also mention that
the corresponding one step mappings $\Phi_n(\eta)=\eta K_{n,\eta}$ and the corresponding
semigroup $\Phi_{p,n}$ satisfy condition $(\Phi_{p,n})$ for some collection of bounded integral operators  $T^{\Phi_{p,n}}_{\eta}$.

\subsection{Description of the main results}\label{descripsec}

The best way to present moderate deviations is to start with the
analysis of the fluctuations of the particle occupation measures.
For mean field particle models, these central limit theorems are
based on a stochastic perturbation interpretation of the local
sampling errors. The random fields associated with these
perturbation models are defined by below.
\begin{defi}
We let $(V^N_n,W^N_n)$ be the sequence of random fields defined by
the pair of  stochastic perturbation formulae:
\begin{equation}\label{defWNn}
\eta^N_n=\eta_{n-1}^NK_{n,\eta_{n-1}^N}+\frac{1}{\sqrt{N}}~V^N_n=\eta_n+\frac{1}{\sqrt{N}}~W^N_n
\end{equation}
where $\eta^N_n=\frac 1N\sum_{j=1}^N \delta_{\xi_n^{(N,j)}}$ is the
empirical distribution of $\xi_n^{N}$.   The sequence $V^N_n$ is
sometimes refereed as the local sampling random field model.
\end{defi}
Notice that the centered random fields $V^N_n$  have conditional
variance functions given by
\begin{equation}\label{covf}
\EE(V_{n}^{N}(f_{n})^2\left|~\xi^{N}_{n-1}~\right.)=\eta_{n-1}^N\left[K_{n,\eta
_{n-1}^{N}}\left((f_n-K_{n,\eta _{n-1}^{N}}(f_{n}))^2\right)\right].
\end{equation}
To analyze the propagation properties of the  sampling errors, {\em
up to a second order remainder measure}, by assumption that
$\Phi_n\in \Upsilon(E_{n-1},E_n)$, we have the first order
decomposition
\begin{equation}\label{decapp}
 \Phi_{n}(\eta)-\Phi_{n}(\mu)\simeq(\eta-\mu) D_{\mu}{\Phi_n}
\end{equation}
with a first order integral operator $D_{\mu}{\Phi_n}$ from
$\Ba(E_n)$ into $\Ba(E_{n-1})$. The precise description of these
regularity properties is provided in section~\ref{regsec}.
\begin{defi}
We let  $(\Da_{p,n})_{0\leq p\leq n}$ be the semigroup $
\Da_{p,n}=\Da_{p+1}\Da_{p+1,n} $, associated with the integral
operator $\Da_n=D_{\eta_{n-1}}{\Phi_n}$. We use the convention
$\Da_{n,n}=Id$, for $p=n$.
\end{defi}

Using the decomposition
\begin{eqnarray}
W^N_n&=&V^N_n+\sqrt{N}\left[\Phi_{n}(\eta_{n-1}^N)-\Phi_n(\eta_{n-1})\right]\nonumber\\
&\simeq& V^N_n+W_{n-1}^ND_{\eta_{n-1}}{\Phi_n} \Longrightarrow
W^N_n\simeq\sum_{p=0}^nV_p^N\Da_{p,n}\label{apporxun}
\end{eqnarray}
 we proved in~\cite{dmrio}
that the sequence of random fields $(V^N_n)_{n\geq 0}$ converges in
law, as $N$ tends to infinity, to the sequence of $n$ independent,
Gaussian and centered random fields $(V_n)_{n\geq 0}$ with a
covariance function with, for any $f,g\in\Ba(E_n)$, the space of the
bounded and measurable real functions on $E_n$ and $n\geq 0$,
\begin{equation}\label{corr1}
\EE(V_n(f)V_n(g))=\eta_{n-1}
K_{n,\eta_{n-1}}([f-K_{n,\eta_{n-1}}(f)][g- K_{n,\eta_{n-1}}(g)]).
\end{equation}
 In addition,  $W^N_n$
converges in law, as  the number of particles $N$ tends to infinity,
to a Gaussian and centered random fields

\begin{equation}\label{W-field}W_n=\sum_{p=0}^nV_p
\Da_{p,n}.\end{equation}

Concentration inequalities associated with these  fluctuations
theorems are developed in some details in a recent article of the
first author with E. Rio~\cite{dmrio}. In this article, we analyze
asymptotic expansions for probabilities of moderate deviations. To
describe with some precision our main results, we recall the
definition of a large deviation principle ({\em abbreviate LDP}) for
random variables.
\begin{defi}
Let $(\alpha(N))_{N\geq 1}$ be a sequence of positive numbers such
that $\lim_{N\rightarrow\infty}\alpha(N)=\infty$. A sequence of
random variables  $\Xa^N$ with values in a topological state space
$(S,\Sa)$ satisfies an LDP with speed $\alpha(N)$ and with good rate
function $I~:~x\in S\mapsto I(x)\in [0,\infty]$ if the pair of
conditions below are satisfied:
\begin{itemize}
\item For every finite constant $a<\infty$, the level sets $\{x\in S~:~I(x)\leq a\}$ are compact sets.
\item For each $A\in\Sa$
$$
-I\left(\stackrel{~o}{A}\right)\leq
\liminf_{N\rightarrow\infty}\frac{1}{\alpha(N)}\log\PP\left(\Xa^N\in
A\right)\leq
\liminf_{N\rightarrow\infty}\frac{1}{\alpha(N)}\log\PP\left(\Xa^N\in
A\right)\leq -I\left(\overline{A}\right)
$$
where, for a subset $B\subset S$, we let $I(B):=\inf_{x\in B}I(x)$.
\end{itemize}
A sequence of random variables $\Ya^N$ is said to satisfy a moderate
deviation principle   ({\em abbreviate MDP})  with good rate
function $I$ and speed $\alpha(N)$ if sequence of random variables
$\Xa^N:={\Ya^N}/{\sqrt{\alpha(N)}} $ satisfies an LDP with speed
$\alpha(N)$ and with good rate function $I$.
\end{defi}

The first steps in the MDP analysis of the pair of random field
sequences $(V^N_n,W^N_n)$  rely on the following pair of asymptotic
Laplace expansions.

\begin{theo}\label{theomdp}
For any nondecreasing function $\alpha(N)$ such that
$\lim_{N\rightarrow\infty}\frac{\alpha(N)}{N}=0$, any $n\geq 0$ and
any collection of functions $f_n\in\Ba(E_n)$, with $n\geq 0$, we
have
\begin{equation}\label{lapV}
\lim_{N\rightarrow\infty}\frac{1}{\alpha(N)}\log{\EE\left(\exp{\left\{\sqrt{\alpha(N)}~\sum_{p=0}^n~
V^N_p(f_p)\right\}}\right)}=\frac{1}{2}\sum_{p=0}^n\EE\left(
V_p\left(f_p\right)^2\right)
\end{equation}
and
\begin{equation}\label{lapW}
\lim_{N\rightarrow\infty}\frac{1}{\alpha(N)}\log{\EE\left(\exp{\left\{\sqrt{\alpha(N)}~W^N_n(f_n)\right\}}\right)}=\frac{1}{2}~
\EE\left( W_n\left(f_n\right)^2\right).
\end{equation}
\end{theo}

The detailed proof of the above theorem is provided in
section~\ref{secondsec} and section~\ref{devsec}. We already mention
that the second expansion  (\ref{lapW}) is a more or less direct
consequence of the first one  (\ref{lapV}) combined with the
 first order decomposition (\ref{apporxun}).

Let us examine some direct consequences of the above theorem.  For
any finite subset
$\Fa_n=\{f^1_n,\ldots,f^{d_n}_n\}\subset\Ba(E_n)^{d_n}$, with
$d_n\geq 1$, we consider the projection mapping defined by
$$
\pi_{\Fa_n}~:~\mu \in \Ma(E_n)\mapsto
\pi_{\Fa_n}(\mu)=\left(\mu(f)\right)_{f\in\Fa_n}\in
\RR^{\Fa_n}\simeq\RR^{d_n}.
$$
By a theorem of J. Gartner and R.S. Ellis, using the asymptotic
Laplace expansion (\ref{lapW})  we prove the following corollary.
\begin{cor}
The random sequence $\pi_{\Fa_n}(W^N_n)$ satisfy an MDP principle in
$\RR^{d_n}$ with speed $\alpha(N)$ and with the good rate function
given for any $v\in\RR^{d_n}$ by the following formula
\begin{equation}\label{leprem}
I^{W_n}_{\Fa_n}(v)=\sup_{u\in\RR^{d_n}}{\left(\langle u,v\rangle
-\frac{1}{2}~\EE\left(
 \left(\sum_{i=1}^d u^i~W_n(f^i_n)\right)^2\right)\right)}\quad\mbox{\rm with}\quad
 \langle u,v\rangle:=\sum_{i=1}^{d_n}u^i v^i.
\end{equation}
If the covariance matrix
$C_{\Fa_n}:=\left(\EE\left(W_n\left(f^i_n\right)W_n\left(f^j_n\right)\right)\right)_{1\leq
i,j\leq d_n}$ is invertible, then the rate function
$I^{W_n}_{\Fa_n}$ takes the form
$$
I^{W_n}_{\Fa_n}(v)=\frac{1}{2} \langle v, C_{\Fa_n}^{-1}v\rangle.
$$
\end{cor}
In much the same way, using the asymptotic Laplace expansion
(\ref{lapV})  we readily prove the following corollary.
\begin{cor}
  The random sequences
$
\left[\pi_{\Fa_0}\left(V^N_0\right),\ldots,\pi_{\Fa_n}\left(V^N_n\right)\right]
$ satisfy a MDP principle in $\RR^{d_0+\ldots+d_n}$ with speed
$\alpha(N)$, with the good rate function given for any
$v=(v_0,\ldots,v_n)\in\RR^{d_0+\ldots+d_n}$ by the following formula
$$
I^{V_{[0,n]}}_{\Fa_{[0,n]}}(v)=\sum_{p=0}^nI^{V_p}_{\Fa_p}(v_p)
$$
with the functions $I^{V_n}_{\Fa_n}$ on $\RR^{d_n}$ defined  as
$I^{W_n}_{\Fa_n}$ by replacing in (\ref{leprem})  the field $W_n$ by
$V_n$.
\end{cor}

In the second part of the article, we strengthen these MDP in two
ways. Firstly, we derive the MDP for the random fields sequences on
the set of measures equipped with the  $\tau$ topology.  Our main
result is the following theorem.

\begin{theo}\label{theopr}
We suppose that  the state spaces $E_n$ are Polish spaces (metric,
complete and separable). In this situation, the sequence of random
fields $ \left(V^N_0,\ldots,V^N_n\right) $ satisfy an MDP in the
product space  $\prod_{p=0}^n\Ma(E_p)$ equipped with the product
$\tau$ topology, with speed $\alpha(N)$ and with the good rate
function $I_{[0,n]}$ given for any $\mu=(\mu_p)_{0\leq p\leq n}\in
\prod_{p=0}^n\Ma(E_p)$ by
$$
I_{[0,n]}(\mu)=\sum_{p=0}^nI_p(\mu_p)$$ with the good rate functions
$I_n$ on $ \Ma(E_n)$ defined for any $ \mu_n\in \Ma(E_n)$ by
\begin{equation}\label{I_n}
I_n(\mu_n)=\sup_{f\in\Ba(E_n)}{\left(\mu_n(f)
-\frac{1}{2}~\eta_{n-1}\left(
 K_{n,\eta_{n-1}}\left[f-K_{n,\eta_{n-1}}(f)\right]^2
 \right)
\right)}.
\end{equation}
In addition, the sequence of random fields
 $W^N_n$ satisfies an MDP  in $\Ma(E_n)$ (equipped with the  $\tau$ topology), with speed $\alpha(N)$
 and with the good rate function
\begin{equation}\label{Jn}
J_n(\nu)=\inf\left\{ \sum_{p=0}^nI_p(\mu_p)~:~\mu~\mbox{\rm
s.t.}~\nu=\sum_{p=0}^n\mu_p\Da_{p,n}
\right\}=\sup_{f\in\Ba(E_n)}{\left(
\nu(f)-\frac{1}{2}\EE\left(W_n(f)^2\right)\right)}.
\end{equation}
\end{theo}

A more explicit description of the rate functions $I_n$ in terms of
integral operators norms on Hilbert spaces can be found in
section~\ref{sectau} (see also section~\ref{varform}, in the
appendix).

Our second main result is a functional moderate deviation for
stochastic processes indexed by a separable collection $\Fa_n$  of
measurable functions $f_n:E_n\rightarrow\RR$ such that $\|f_n\|\leq
1$. We let $l_{\infty}(\Fa_n)$ be the space  of
all bounded real functions $F_n~:~f\in\Fa_n\mapsto F_n(f_n)\in\RR$
on $\Fa_n$ with the sup norm
$\|F_n\|_{\Fa_n}=\sup_{f_n\in\Fa_n}|F_n(f_n)|$. Notice that this
vector space is a non separable Banach space if the set of functions
$\Fa_n$ is infinite. To measure the size of a given class $\Fa_n$,
one considers the covering numbers $N(\e,\Fa_n,L_p(\mu))$
\index{Covering numbers}defined as the minimal number of
$L_p(\mu)$-balls of radius $\e >0$ needed to cover $\Fa_n$. By
$\Na(\e,\Fa_n)$, $\e >0$, and by $I(\Fa_n)$ we denote the uniform
covering numbers and entropy integral \index{Entropy integral}given
by
$$
\Na(\e,\Fa_n)=\!\sup_{\eta\in\Pa(E_n)} \bigl
\{\Na(\e,\Fa_n,\LL_2(\eta))\bigr \} \quad\mbox{\rm and}\quad
\Ia(\Fa_n)=\!\int_{0}^2\!\sqrt{\log{ \Na(\e,\Fa_n) }}\,d\e.
$$
We further assume that
\begin{description}
\item{\bf (A1)} $\Na(\e,\Fa_n)<\infty$ for any $\epsilon>0$, and $
\Ia(\Fa_n)<\infty$.
\end{description}
This condition implies that the set $\Fa_n$ is totally bounded in $L_2(\eta)$,
for any distribution $\eta$ on $E_n$. Various examples of classes of
functions with finite covering and entropy integral are given in the
book of Van der Vaart and Wellner~\cite{wellner} (see for instance
p. 86, p. 135, and exercise 4 on p.150).

For any $\delta>0$, we also set
$$
\Fa_n(\delta):=\left\{h=(f-g)~:~(f,g)\in \Fa_n~\mbox{\rm s.t.}~\eta_n(h^2)^{1/2}\leq \delta\right\}
$$

\begin{description}
\item{\bf (A2)}
There exists some
 separable collection $\Fa_{n}^{\prime}$  of
measurable functions $f_n$ on $E_{n}$, s.t. $\|f_n\|\leq
1$, s.t. $I(\Fa^{\prime}_{n})<c_0(n)~I(\Fa_{n+1})$, and such that for any
probability measure $\mu$, any $\delta>0$, we have
$$
\left\|\Phi_{n+1}(\mu)-\Phi_{n+1}(\eta_{n})\right\|_{\Fa_{n+1}(\delta)}
\leq c_2(n)~\left\|\mu-\eta_{n}\right\|_{\Fa^{\prime}_{n}(c_1(n)\delta)}
$$
for some finite constant $c_i(n)<\infty$, $i=0,1,2$, whose values only depend on the
mapping $\Phi_{n+1}$, and on the measure $\eta_{n}$.
\end{description}
We illustrate this regularity condition in the context of the Feynman-Kac models
presented in (\ref{phiFK}). Using (\ref{ref-lipsch}), we find that
$$
\left|\left[\Phi_{n+1}\left(\eta\right)-\Phi_{n+1}\left(\eta_n\right)\right](h)\right|\leq g_n~\left|\left(\mu-\eta_n\right)
\left(
\frac{G_n}{\eta_n(G_n)}~\left(M_{n+1}(h)-\eta_{n+1}(h)\right)
\right)
\right|
$$
where $g_n=\sup_{x,y}G_n(x)/G_n(y)$ and
\begin{eqnarray*}
\eta_n\left(\left(
\frac{G_n}{\eta_n(G_n)}~\left(M_{n+1}(h)-\eta_{n+1}(h)\right)
\right)^2\right)&\leq& g_n~\eta_n\left(
\frac{G_n}{\eta_n(G_n)}~\left(\left(M_{n+1}(h)-\eta_{n+1}(h)\right)
\right)^2\right)\\
&\leq &g_n~\eta_{n+1}(h^2)
\end{eqnarray*}
Using elementary manipulations, we show that ${\bf (A2)}$ is met with the constants
$c_1(n)=1/(2\sqrt{g_n})\leq 1$,  $c_2(n)=2g_n^2$ and the class of functions
$$
\Fa^{\prime}_{n}=\left\{
\frac{1}{2g_n}~\frac{G_n}{\eta_n(G_n)}~\left(M_{n+1}(f)-\eta_{n+1}(f)\right)~:~f\in \Fa_{n+1}
\right\}
$$
Using lemma 2.3 in~\cite{dm-ledoux}, we also prove that $I(\Fa^{\prime}_{n})<c_0(n)~I(\Fa_{n+1})$ for some
finite constant whose values only depends on $g_n$.

For any finite subset
$\Ga_n\subset\Fa_n$, we let
$$
\pi_{\Fa_n,\Ga_n}~:~v \in l_{\infty}(\Fa_n)\mapsto
\pi_{\Fa_n,\Ga_n}(v)=\left(v(g)\right)_{g\in\Ga_n}\in
l_{\infty}(\Ga_n)=\RR^{\Ga_n}
$$
be the restriction mapping defined by
$\pi_{\Fa_n,\Ga_n}(\nu)(g)=v(g)$, for any $g_n\in\Ga_n$. The MDP of
the stochastic processes $W^N_n$ on  $\La_{\infty}(\Fa_n)$ are
described below.
\begin{theo}\label{theoprpp} Assume that the class of observables $\Fa_n$
satisfies (A1), and condition (A2) is met. The sequence of stochastic processes $W^N_n$ satisfy
the large deviation principle in
 $\La_{\infty}(\Fa_n)$ with the good rate function $I_{\Fa_n}^{W_n}$ given below
\begin{eqnarray*}
v\in \La_{\infty}(\Fa_n) ~~~~~~ I_{\Fa_n}^{W_n}(v)&=\sup\left\{
I^{W_n}_{\Ga_n}(\pi_{\Fa_n,\Ga_n}(v))~:~\Ga_n\subset\Fa_n~,~
\mbox{with}~\Ga_n~\mbox{ finite} \right\}\\
&=\inf\left\{J_n(\nu)|\nu\in M_0(E_n),\nu(f)=v(f),\forall f\in
\Fa_n\right\}.
 \end{eqnarray*}
 where $J_n$ is given in (\ref{Jn}).
\end{theo}

For finite sets $\Fa_n$, the above theorem clearly reduces to the
MDP presented in (\ref{leprem}). Also observe that
 the
$\tau$-topology  on $\Ma(E_n)$ is sometimes finer than the topology
associated with the seminorm $\|\mu-\eta\|_{\Fa_n}$ induced by
$\Fa_n$. For instance, when $E=\RR^d$ and
$\Fa=\{1_{(-\infty,x]}\;;\;x\in \RR^d\}$, the topology induced by
the supremum distance
$$
\|\mu-\eta\|_{\Fa}=\sup_{x\in\RR^d}{|\mu((-\infty,x])-\eta((-\infty,x])|}
$$
is strictly coarser than the $\tau$-topology. In this situation,
Theorem~\ref{theoprpp} is a direct consequence of
Theorem~\ref{theopr}. In more general situations, by \cite{liming2}
or a theorem~of M. A. Arcones (see for instance~theorem 3.2
in~\cite{arcones}), the MDP for stochastic processes $W^N_n$ in
$\La_{\infty}(\Fa_n)$ is deduced from the MDP of the finite
marginals $\pi_{\Fa_n,\Ga_n}(W^N_n)$ plus the following exponential
asymptotic equicontinuity condition:
$$
\forall y>0,\qquad \lim_{\delta\rightarrow
0}\limsup_{N\rightarrow\infty}\frac{1}{\alpha(N)}\log{
\PP\left(\frac{1}{\sqrt{\alpha(N)}}~\left\|W_n^N\right\|_{\Fa_n(\delta)}>y\right)}=-\infty
$$
with the collection of functions
$$
\Fa_n(\delta):=\{h_n~:~h_n=(f_n-g_n)~\mbox{\rm
with}~(f_n,g_n)\in\Fa_n^2~\mbox{\rm and}~\eta_n(h_n^2)\leq \delta\}.
$$

\section{Asymptotic Laplace expansions}

\subsection{Some preliminary results}

 \begin{lem}
 For any $0\leq p\leq n$, we have
 $\Phi_{p,n}\in  \Upsilon(E_{p},E_n)$ with  the first order
decomposition type formula
\begin{equation}\label{decompPhi}
 \Phi_{p,n}(\eta)-\Phi_{p,n}(\mu)=[\eta-\mu]D_{\mu}\Phi_{p,n}+\Ra^{\Phi_{p,n}}(\eta,\mu)
\end{equation}
 for some collection of bounded integral kernels $D_{\mu}\Phi_{p,n}$  from $E_p$ into $E_n$ and
 some second order remainder signed measures $\Ra^{\Phi_{p,n}}(\eta,\mu)$.
 In addition, for any $N\geq 1$, we have the first order decomposition
\begin{equation}\label{decompV}
W^N_n=\sum_{p=0}^nV^N_p\Da_{p,n}+\frac{1}{\sqrt{N}}~\Ra^N_n\quad\mbox{with}\quad
\Ra^N_n:=
N~\sum_{p=0}^{n-1}R_{p+1}^{\Phi_{p+1}}\left(\eta^N_{p},\eta_p\right)D_{p+1,n}
\end{equation}
and the semigroup $(\Da_{p,n})_{0\leq p\leq n}$ introduced in
(\ref{apporxun}).
\end{lem}

\begin{lem}\label{lem1b}
 For every $f\in \mbox{\rm Osc}_1(E_n)$, $N\geq 1$ and any $n\geq 0$ and $m\geq 1$, we have
 the  $\LL_{m}$ estimates:
 \begin{equation}\label{eqlm}
 \EE\left(
 \left|V_{n}^{N}(f_{n})\right|^{m}
 \left|
 \xi^{(N)}_{n-1}
 \right. \right)^{\frac{1}{m}}\leq
b(m)\quad\mbox{ and}\quad\sqrt{N}~
\EE\left(\left|\left[\eta^N_n-\eta_n\right](f_n)\right|^{m}
\right)^{\frac{1}{m}}\leq b(m)~\sum_{p=0}^n\delta(T^{\Phi_{p,n}})
\end{equation}
 as well as the bias estimate
 \begin{equation}\label{biasprop}
N~\left|\EE\left(\eta^N_n(f_n)\right)-\eta_n(f_n)\right|\leq
\sum_{p=0}^n  \delta(R^{\Phi_{p,n}}).
\end{equation}
\end{lem}

A detailed proof of (\ref{decompPhi})
 can be found in~\cite{dmrio}.
Formula (\ref{decompV}) is a direct consequence of the following
inductive decomposition
$$
W^N_n=V^N_n+W^N_{n-1}\Da_n+\sqrt{N}~R^{\Phi_n}\left(
\eta^N_{n-1},\eta_{n-1} \right).
$$
The proof of lemma~\ref{lem1b} is postponed to
section~\ref{seclem1b} in the appendix.

\subsection{Second order remainder measures}\label{secondsec}

This section is mainly concerned with non asymptotic Laplace estimates of the second
order remainder measures  introduced in lemma~\ref{decompV}, namely
$$\Ra^N_n:=
\sqrt{N}\left[ W^N_n-\sum_{p=0}^nV^N_p\Da_{p,n} \right].
$$

\begin{prop}\label{propeqq}
For every $f\in \mbox{\rm Osc}_1(E_n)$, $N\geq 1$, $n\geq 0$,  we have the Laplace
estimates :
\begin{equation}\label{refequiv}
\forall
t\in\left[0,1/(2r(n))\right)\qquad\EE\left(\exp{\left(t~\sqrt{N}\left|\Ra^N_n(f_n)\right|\right]}\right)\leq
\frac{1}{\sqrt{1-2r(n)t}}
\end{equation}
with some finite constant
$
 r(n)\leq \sum_{p=0}^{n-1}\beta(\Da_{p+1,n})~\left(\sum_{q=0}^p\delta(T^{\Phi_{q,p}})\right)^2~\delta\left(R^{\Phi_{p+1}}\right)
$.

\end{prop}

\proof
By (\ref{decompV}), we have that
$$
\left|\Ra^N_n(f_n)\right|\leq \sum_{p=0}^{n-1}
\int~\left|\left(V^N_p\right)^{\otimes
2}(g)\right|~R^{\Phi_{p+1}}_{\eta_{p}}(f,dg).
$$
Combining (\ref{eqlm}) with the generalized Minkowski inequality this implies that
$$
\left(\EE\left|\Ra^N_n(f_n)\right|^m\right)^{1/m}\leq b(2m)^2~r(n).
$$
We end the proof of the proposition recalling that for Gaussian centered random variable with $\EE(X^2)=1$
we have that $b(2m)^{2m}=\EE(X^{2m})$ and for any $t\in [0,1/2[$
$$
\EE(\exp{\left\{tX^2\right\}})=\sum_{m\geq 0}\frac{t^m}{m!}
~b(2m)^{2m}= {1}/{\sqrt{1-2t}}.
$$
\cqfd

\begin{cor}
For every $f\in \mbox{\rm Osc}_1(E_n)$, $N\geq 1$, $n\geq 0$, and for every  $\epsilon>0$, we have
$$
\PP\left(\left|\Ra^N_n(f_n)\right|\geq \epsilon+\frac{r(n)}{\sqrt{N}}\right)\leq 2
 e^{-\frac{\epsilon\sqrt{N}}{2r(n)}\left\{1
- \delta_n(\epsilon, N)
 \right\}}
~~ \mbox{where }~~ \delta_n(\epsilon,
N)=\frac{r(n)}{\epsilon\sqrt{N}}\log{\left(1+\frac{\epsilon\sqrt{N}}{r(n)}\right)}.
$$
In particular,  for any nondecreasing function $\alpha(N)$ such that
$\lim_{N\rightarrow\infty}\frac{\alpha(N)}{N}=0$, we have
\begin{equation}\label{asympeq}
 \lim_{N\rightarrow\infty}\frac{1}{\alpha(N)}\log{\PP\left(
\left|\Ra^N_n(f_n)\right|
\geq \epsilon \sqrt{\alpha(N)}
 \right) }=-\infty.
\end{equation}
In other words,
the random fields $\frac{1}{\sqrt{\alpha(N)}}W^N_n$ and $\frac{1}{\sqrt{\alpha(N)}}\sum_{p=0}^n ~V_p^N \Da_{p,n}$
are $\alpha(N)$-exponentially equivalent.
\end{cor}
\proof
Using the fact that
$$
\log{\EE\left(e^{t\left[\Ra^N_n(f_n)-r(n)\right]}\right)}\leq
-r(n)t-\frac{1}{2}\log{\left(1-2r(n)t\right)},
$$
we readily find that
$$
\PP\left(\Ra^N_n(f_n)\geq \epsilon+r(n)\right)\leq
\exp{\left(-\sup_{t\leq
1/2}\left\{\frac{\epsilon}{r(n)}~t+t+\frac{1}{2}\log{(1-2t)}\right\}\right)}.
$$
Choosing $t=\frac{1}{2}\left(1-\frac{1}{1+\epsilon}\right)$, we find that
$$
\PP\left(\Ra^N_n(f_n)\geq \epsilon+r(n)\right)\leq
 \exp{\left(-\frac{\epsilon}{2r(n)}\left\{1-\frac{r(n)}{\epsilon}\log{\left(1+\frac{\epsilon}{r(n)}\right)}\right\}\right)}
$$
which ends the proof of the corollary.\cqfd

We end this section with a technical  transfer lemma of Laplace
asymptotic expansions for  arbitrary stochastic processes. The proof
is elementary, so omitted. 

\begin{lem}\label{transferlemma}
Let $(X_N)$, $(Y_N)$ two sequences of random valuables such that for
any $\lambda\geq0$,
$$\lim_{N\rightarrow\infty}\frac{1}{\alpha(N)}\log\EE\left(e^{\lambda\alpha(N)X_N}\right)=\Lambda(\lambda)\quad and\quad
\lim_{N\rightarrow\infty}\frac{1}{\alpha(N)}\log\EE\left(e^{\lambda\alpha(N)|X_N-Y_N|}\right)=0
$$
for some sequence $\alpha(N)$ increasing to infinite and some finite
logarithmic moment generating function $\Lambda(\lambda)$. Then for
all $\lambda\ge0$, we have
$$\lim_{N\rightarrow\infty}\frac{1}{\alpha(N)}\log\EE(e^{\lambda\alpha(N) Y_N})=\Lambda(\lambda).$$
\end{lem}

\subsection{Asymptotic Laplace transform estimates}\label{devsec}

This section is mainly concerned with the proof of
Theorem~\ref{theomdp}. The fluctuation properties  of the first
order random field sequence $\sum_{p=0}^nV^N_p\Da_{p,n}$ is encoded
in the pair of martingale sequences defined below.
\begin{defi}
We associate with collection of functions $f=(f_n)_{n\geq
0}\in\prod_{n\geq 0}\Ba(E_n)$, the pair of
$\sigma\left(\xi^{(N)}_0,\ldots,\xi^{(N)}_n\right)$-martingale
sequences given below
$$
M^{(N)}_n(f)=\sum_{p=0}^nV^N_p(f_p)\quad \mbox{and}\quad
E^{(N)}_n(f):=\frac{1}{\Za^{(N)}_n(f)}~\exp{\left\{\sqrt{\alpha(N)}~ M^{(N)}_n(f)\right\}}
$$
with the stochastic product
$$
\Za^{(N)}_n(f):=\prod_{p=1}^n
\EE\left(\exp{\left\{\sqrt{\alpha(N)}~V^N_p(f_p)\right\}}~|~\xi^{(N)}_{p-1}\right).
$$

\end{defi}
For every $N\geq 1$, we notice that the angle bracket of $M^{(N)}_n(f)$ is given by
$$
\langle M^{(N)}(f)\rangle_n=\sum_{p=0}^n \Delta_p\langle M^{(N)}(f)\rangle
$$
with the random increments
$$
\Delta_n\langle M^{(N)}(f)\rangle
:=\eta^N_{n-1}\left(K_{n,\eta^N_{n-1}}\left[\left(f_n-K_{n,\eta^N_{n-1}}(f_n)\right)^2\right]\right).
$$
We know that the sequence of martingales $ M^{(N)}_n(f) $ converges
in law, as $N$ tends to infinity, to the Gaussian martingale
$$
M_n(f)=\sum_{p=0}^n ~V_p(f_p) \quad\mbox{\rm with}\quad \langle
M(f)\rangle_n =\sum_{p=1}^n
\eta_{p-1}\left(K_{p,\eta_{p-1}}\left[\left(f_p-K_{p,\eta_{p-1}}(f_p)\right)^2\right]\right).
$$
The main object of this subsection is to prove that

\begin{equation}\label{laplaceest}
\lim_{N\rightarrow\infty}\frac{1}{\alpha(N)}\log{\EE\left(e^{\alpha(N)\left(\frac{1}{\sqrt{\alpha(N)}}
~M^{(N)}_n(f)\right)}\right)}=\frac{1}{2}~\langle
M(f)\rangle_n.
\end{equation}

Notice that the above asymptotic Laplace expansion is equivalent to (\ref{lapV}).

The next technical lemma is pivotal.

\begin{lem}\label{lem3}
There exist a pair of functions $(\tau^{(N)}_{j,n}(f))_{j=1,2}$ that
converge to $0$ as $N$ tends to $\infty$, such that
$$
e^{\sqrt{\alpha(N)} M^{(N)}_n(f)-\frac{\alpha(N)}{2}~\langle M^{(N)}(f)\rangle_n}
~\leq~ E^{(N)}_n(f)~~e^{\tau^{(N)}_{2,n}(f)~\frac{\alpha(N)}{2}~\langle M^{(N)}(f)\rangle_n}
$$
and
$$
E^{(N)}_n(f)~~e^{-\tau^{(N)}_{1,n}(f)~\frac{\alpha(N)}{2}~\langle
M^{(N)}(f)\rangle_n}~\leq~
e^{\sqrt{\alpha(N)}M^{(N)}_n(f)-\frac{\alpha(N)}{2}~\langle
M^{(N)}(f)\rangle_n}.
$$
\end{lem}

The  proof of lemma~\ref{lem3} is rather technical, thus we postpone
it to section~\ref{prooflem3} in the appendix.

\begin{prop}\label{prop1}
$$
\lim_{N\rightarrow\infty}\frac{1}{\alpha(N)}\log{\EE\left(e^{\alpha(N)\left(\frac{1}{\sqrt{\alpha(N)}}
~M^{(N)}_n(f)-\frac{1}{2}~\langle M^{(N)}(f)\rangle_n\right)}\right)}=0
$$
and
 \begin{equation}\label{mmN}
\EE\left(e^{t  \sqrt{N}\left|\langle M^{(N)}(f)\rangle_n-\langle
M(f)\rangle_n\right|}\right)\leq \left(1+t\overline{c}_n \right)~
e^{{(\overline{c}_nt)^2}/{2} }.
\end{equation}
In the above display, $\overline{c}_n$ stands for some finite
constant $\overline{c}_n:= \sum_{p=0}^nc(p)$ with
$$
 c(p):=2\left\{1+
\delta\left(T^{\Phi_p}\right)+ \delta\left(T^{K_p}\right)
\right\}~\sum_{0\leq q<p}\delta(T^{\Phi_{q,p-1}}).~$$
\end{prop}

Before getting into the details of the proof of the above proposition, it is convenient to make a couple
of comments. Firstly, replacing in (\ref{mmN}) the parameter $t$ by $\frac{\alpha(N)}{\sqrt{N}}t$ we find that
 \begin{eqnarray*}
\EE\left(e^{t  \alpha(N)\left|\langle M^{(N)}(f)\rangle_n-\langle
M(f)\rangle_n\right|}\right)&\leq
&\left(1+\frac{t\alpha(N)}{\sqrt{N}}~\overline{c}_n \right)~
\exp{\left\{\frac{t^2\alpha(N)^2}{2N}~\overline{c}_n^2 \right\}}.
\end{eqnarray*}
from which we conclude that
$$
\forall t\geq 0, \qquad
\limsup_{N\rightarrow\infty}\frac{1}{\alpha(N)}{\EE\left(e^{\alpha(N)~t\left|\langle
M^{(N)}(f)\rangle_n-\langle M(f)\rangle_n\right|}\right)}=0.
$$
Also observe that the stochastic processes
\begin{eqnarray*}
A^N_n(f)&=&\frac{1}{\sqrt{\alpha(N)}}~M^{(N)}_n(f)-\frac{1}{2}~\langle M^{(N)}(f)\rangle_n\\
B^N_n(f)&=&\frac{1}{\sqrt{\alpha(N)}}~M^{(N)}_n(f)-\frac{1}{2}~\langle M(f)\rangle_n
\end{eqnarray*}
on the set of sequence $f=(f_p)_{0\leq p\leq n}\in\prod_{p=0}^n\Ba(E_p)$, have the following scaling properties
$$
\left|A^N_n(f)-\epsilon^{-1}~A^N_n(\epsilon f)\right|=\frac{1}{2}~\langle M^{(N)}(f)\rangle_n~(1-\epsilon)
\leq \frac{1}{2}~(1-\epsilon)~\sum_{p=0}^n\mbox{\rm osc}(f_p)^2
$$
and
$$
\left|B^N_n(f)-\epsilon^{-1}~B^N_n(\epsilon f)\right|=\frac{1}{2}~\langle M(f)\rangle_n~(1-\epsilon)
$$
for any $\epsilon\in [0,1]$. In the above display, $\epsilon f$ stands for the
sequence of functions $(\epsilon f_p)_{0\leq p\leq n}$.
Therefore the asymptotic Laplace expansion (\ref{laplaceest}) is  a direct consequence of the transfer lemma~\ref{transferlemma}.

Now, we come to

{\bf Proof of proposition \ref{prop1}.}

Since we have
$
\langle M^{(N)}(f)\rangle_n\leq \sigma_n^2(f):=\sum_{p=0}^n\mbox{\rm osc}(f_p)^2
$,
using lemma~\ref{lem3} we readily prove that
$$
-\tau^{(N)}_{1,n}(f)~\frac{1}{2}~\sigma_n^2(f)\leq
\frac{1}{\alpha(N)}\log{\EE\left(e^{\alpha(N)\left(\frac{1}{\sqrt{\alpha(N)}}~M^{(N)}_n(f)-\frac{1}{2}~\langle
M^{(N)}(f)\rangle_n\right)}\right)}\leq
\tau^{(N)}_{2,n}(f)~\frac{1}{2}~\sigma_n^2(f).
$$
This ends the proof of the first assertion. Now, we come to the proof of (\ref{mmN}).
For every $n\geq 1$, $\eta\in\Pa(E_{n-1})$ and $f_n\in \Ba(E_n)$ we set
$$
\Sigma_n(\eta,f_n)
:=\eta\left(K_{n,\eta}\left[\left(f_n-K_{n,\eta}(f_n)\right)^2\right]\right).
$$
For $n=0$, we set $\Sigma_0(\eta,f_0)=\eta([f_0-\eta(f_0)]^2)$.
Firstly, we observe that
\begin{eqnarray*}
\Sigma_n(\eta,f_n)-\Sigma_n(\mu,f_n)
&=&\left[\Phi_n(\eta)-\Phi_n(\mu)\right]\left(f_n^2\right)+\mu\left(K_{n,\mu}(f_n)^2\right)-
\eta\left(K_{n,\eta}(f_n)^2\right)\\
&=&\left[\Phi_n(\eta)-\Phi_n(\mu)\right]\left(f_n^2\right)+\left[\mu-\eta\right]\left(K_{n,\eta}(f_n)^2\right)\\
&&\hskip4cm + \mu\left(K_{n,\mu}(f_n)^2-K_{n,\mu}(f_n)^2\right).
\end{eqnarray*}
This implies that
$$
\begin{array}{l}
\vert\Sigma_n(\eta,f_n)-\Sigma_n(\mu,f_n)\vert\\
\\
\leq
\vert\left[\Phi_n(\eta)-\Phi_n(\mu)\right]
\left(f_n^2\right)\vert+\vert\left[\mu-\eta\right]\left(K_{n,\eta}(f_n)^2\right)\vert+2\Vert K_{n,\mu}(f_n)-K_{n,\eta}(f_n)\Vert
\end{array}
$$
and therefore
 $$
 \begin{array}{l}
 \left(\EE\vert
 \Sigma_n(\eta^N_{n-1},f_n)-\Sigma_n(\eta_{n-1},f_n)
 \vert^m\right)^{\frac{1}{m}}\\
 \\
\leq  \displaystyle\int~
\left(\EE\vert(\eta^N_{n-1}-\eta_{n-1})(g)\vert^m\right)^{\frac{1}{m}}~T^{\Phi_n}_{\eta_{n-1}}(f_n^2,dg)
+
\left(\EE\vert(\eta^N_{n-1}-\eta_{n-1})\left(K_{n,\eta_{n-1}}(f_n)^2\right)\vert^m\right)^{\frac{1}{m}}\\
\\
\hskip5cm +2\displaystyle\int~
\EE\left(\vert(\eta^N_{n-1}-\eta_{n-1})(g)\vert^m\right)^{\frac{1}{m}}~T^{K_n}_{\eta_{n-1}}(f_n,dg).
 \end{array}
 $$
 Using (\ref{eqlm}), we have the upper bound
 $$
 \sqrt{N}~\EE\left(\vert
 \Sigma_n(\eta^N_{n-1},f_n)-\Sigma_n(\eta_{n-1},f_n)
 \vert^m\right)^{\frac{1}{m}}
\leq b(m)~c(n).$$

 One concludes that
 $$
 \sqrt{N}~\EE\left(\vert
\langle M^{(N)}(f)\rangle_n-\langle M(f)\rangle_n
 \vert^m\right)^{\frac{1}{m}}\leq b(m)~\overline{c}_n .$$
The $\LL_m$-inequalities stated above
clearly imply that for any  $t>0$
$$
\begin{array}{l}
\EE\left(\exp{\left\{t \sqrt{N}\left| \langle M^{(N)}(f)\rangle_n-\langle M(f)\rangle_n\right|\right\}}\right)\\
\\
=\sum_{m\geq 0}\frac{t^{2m}}{(2m)!}\EE\left(\left(
 \langle M^{(N)}(f)\rangle_n-\langle M(f)\rangle_n\right)^{2m}\right)\\
 \\
 \hskip4cm+\sum_{m\geq 0}\frac{t^{2m+1}}{(2m+1)!}\EE\left(\left| \langle M^{(N)}(f)\rangle_n-\langle M(f)\rangle_n\right|^{2m+1}\right)\\
\\
\leq \sum_{m\geq 0}\frac{1}{m!}~
\left(\frac{t^2\overline{c}_n^2}{2}\right)^m+(t\overline{c}_n)
\sum_{m\geq 0}\frac{1}{m!}~
\left(\frac{t^2\overline{c}_n^2}{2}\right)^{m}.
\end{array}
$$
where (\ref{mmN}) follows. \cqfd

\subsection{Proof of Theorem \ref{theomdp}}

{\it Proof of (\ref{lapV}).} This is done in Subsection 3.3.

\medskip\noindent
{\it Proof of (\ref{lapV})$\implies$ (\ref{lapW}). } Note that if a
final time horizon $n$ is fixed then we have for any function
$f_n\in\Ba(E_n)$
$$
\left(\forall 0\leq p\leq n\quad
f_p=\Da_{p,n}(f_n)\right)\Longrightarrow
\sum_{p=0}^nV^N_p(f_p)=\sum_{p=0}^n ~V_p^N \Da_{p,n}(f_n).
$$
Let $(A^N_n,B^N_n)$ the pair of random fields defined below:
$$A^N_n=\frac{1}{\sqrt{\alpha(N)}}\sum_{p=0}^n ~V_p^N \Da_{p,n}\quad \mbox{\rm and}
\quad B^N_n=\frac{1}{\sqrt{\alpha(N)}}~W^N_n.
$$
By (\ref{lapV}), we have
$$
\lim_{N\rightarrow\infty}\frac{1}{\alpha(N)}\log{\EE\left(e^{\alpha(N)~A^N_n(f_n)}\right)}=
A_n(f_n):=\EE\left(
\frac{1}{2}\sum_{p=0}^nV_p(\Da_{p,n}(f_n))^2\right)
$$
and by (\ref{refequiv})
$$
\forall
t\in\left[0,N/(2\alpha(N)r(n))\right[,\qquad\EE\left(e^{t~\alpha(N)\left|\left[
B^N_n-A^N_n \right](f_n)\right|}\right)\leq
\left(1-\frac{\alpha(N)2r(n)t}{N}\right)^{-\frac{1}{2}}.
$$
This yields that
$$
\forall t>0\qquad
\lim_{N\rightarrow\infty}\frac{1}{\alpha(N)}\log{\EE\left(e^{~t\alpha(N)\left|\left[
B^N_n-A^N_n \right](f_n)\right|}\right)}=0
$$
where (\ref{lapW}) follows by the transfert
lemma~\ref{transferlemma}.

\section{Moderate deviations in $\tau$-topology}

We further require that the state spaces $E_n$ are Polish spaces.
The $\tau$-topology on $\Ma(E_n)$ is the coarsest topology that
makes the maps $\mu\in\Ma(E_n)\mapsto \mu(A)$ continuous, for any
measurable set $A\in\Ea$.

\subsection{A deviation theorem for the local sampling random fields}\label{sectau}
The main object of this section is to prove the following theorem.
\begin{theo}\label{theopr2}
The sequence of random fields $ V^N_n $ satisfy  an MDP in
$\Ma(E_n)$ equipped with the $\tau$ topology, with speed $\alpha(N)$
and with the good rate function
\begin{equation}\label{casun}
I_n(\mu)=\sup_{f\in\Ba(E_n)}{\left(\mu(f) -\frac{1}{2}~\eta_{n-1}\left(
 K_{n,\eta_{n-1}}\left[f-K_{n,\eta_{n-1}}(f)\right]^2
 \right)
\right)}.
\end{equation}
In addition, for any $n\geq 0$, the sequence of random fields
$
V^N_{[0,n]}:=(V^N_0,\ldots,V^N_n)
$
satisfy  an MDP in the product space $\prod_{p=0}^n\Ma(E_n)$, with speed $\alpha(N)$ and with the good rate function
$$
I_{[0,n]}(\mu_0,\ldots,\mu_n)=\sum_{p=0}^nI_p(\mu_p).
$$
\end{theo}
Before entering into the proof of this theorem, we provide a more explicit representation
of the rate functions $I_n$. Firstly, assume that the McKean transitions $K_{n,\eta}$ are given by
$K_{n,\eta}(x,dy)=\Phi_{n}(\eta)(dy)$. In this situation, we have
$$
K_{n,\eta_{n-1}}(x,dy)=\eta_n(dy)\Rightarrow
I_n(\mu)=\sup_{f\in\Ba(E_n)}{\left(\mu(f) -\frac{1}{2}~\eta_{n}\left(
\left[f-\eta_n(f)\right]^2
 \right)
\right)}.
$$
The variational formula given above coincides with the one of the rate function
of  the MDP associated with independent and identically distributed  random sequences.
In this case, we have that
$$
I_n(\mu)=
\left\|\frac{d\mu}{d\eta_n}\right\|_{\LL_2(\eta_n)}^2\quad \mbox{\rm if}\quad
 \mu\ll\eta_n\quad\mbox{\rm with}\quad \frac{d\mu}{d\eta_n}\in \LL_2(\eta_n)\quad\mbox{and} \quad \mu(E)=0
$$
and $I_n(\mu)=\infty$, otherwise. A proof of this assertion is provided in section~\ref{varform}, in the appendix.
In more general situations, we need to work a little harder.
 Let  $K^{\star}_{n,\eta_{n-1}}$ be the adjoint  operator of $K_{n,\eta_{n-1}}$ from
$\LL_2(\eta_{n-1})$ into $ \LL_2(\eta_n)$ given by
$$
\forall (f,g)\in  \LL_2(\eta_{n})\times \LL_2(\eta_{n-1}),\quad
\eta_n\left(f K^{\star}_{n,\eta_{n-1}}(g)\right)=\eta_{n-1}
(K_{n,\eta_{n-1}}(f)~g).
$$
We will prove in section 7.4. the following explicit expression:
\begin{equation}
I_n(\mu)=\frac{1}{2} \sum_{m\geq 0} \eta_n\left[ \left(h_\mu\right)
\left(K^{\star}_{n,\eta_{n-1}}K_{n,\eta_{n-1}}\right)^m
 \left(h_\mu \right)\right], \ \mu\in {\mathcal M}_0(E_n),
 \mu\ll\eta_n, h_\mu =\frac{d\mu}{d\eta_n}\in \LL_2(\eta_n)
\end{equation}
and $I_n(\mu)=+\infty$ otherwise.

\subsection{Moderate deviations for projective limits}

The proof of the theorem is based on a projective limit
interpretation of the strong topology on the set of finite and
signed measures. We begin by first introducing several definitions.

\begin{defi}
We let $\Ua(E_n)$ the set of finite partitions $U_{n}=(U^i_n)_{1\leq
i\leq d}\in \Ea^d_n$ of the set $E_n$, with $d\geq 1$. We let
$\sigma(U_{n})$ be the $\sigma$-field generated by $U_{n}$.
 We also let $$\pi_{U_{n}}~:~\mu\in\Ma(E)\mapsto\pi_{U_{n}}(\mu)\in\Ma(E_n,\sigma(U_{n}))$$
be the restriction of the measure $\mu$ to the sigma-field
$\sigma(U_{n})$.
\end{defi}
Notice that $\Ma(E_n,\sigma(U_{n}))$ can be identified with
$\RR^{U_{n}}\simeq\RR^d$. Furthermore, the $\sigma$-algebra and the
$\tau$-topology induced on $\Ma(E_n,\sigma(U_{n}))$ by the
restriction mapping $\pi_{U_{n}}$ coincide with the natural topology
and  the Borel sigma-field on $\RR^d$.

\begin{defi}
 We say that a partition
$U^{\prime}_n$ is finer than $U_{n}$, and we write $U^{\prime}_n\geq
U_n$, as soon as we have $\sigma(U^{\prime}_n)\supset \sigma(U_n)$.
We also let
$\pi_{U_n^{\prime},U_n}~:~\mu\in\Ma(E,\sigma(U^{\prime}_n))\mapsto\pi_{U^{\prime}_n,U_n}(\mu)\in\Ma(E_n,\sigma(U_n))$
be the restriction of the measure $\mu$ on $\sigma(U^{\prime}_n)$ to
the sigma-field $\sigma(U_n)$. The set
$\left(\Ma(E_n,\sigma(U_n)),\pi_{U_n^{\prime},U_n}\right)_{U^{\prime}_n\geq
U_n}$ forms a projective inverse spectrum of $\Ua(E_n)$. We let
$\lim_{\Ua_n}\Ma_n$ be the projective limit space of the spectrum
$$
\lim_{\Ua_n}\Ma_n:=\left\{\mu\in\prod_{U_n\in\Ua_n}\Ma(E_n,\sigma(U_n))~:~\forall
U^{\prime}_n\geq U_n\quad \pi_{U_n}(\mu)=
\pi_{U_n^{\prime},U_n}(\pi_{U_n^{\prime}}(\mu) )\right\}.
$$
\end{defi}

\begin{defi}
We let $\Mb (E_n)$ be the set of finite  additive set functions from
$\Ea_n$ into $\RR_+$, equipped with the $\tau_1$-topology of setwise
convergence. More precisely, a sequence $\mu_k\in \Mb(E_n)$
$\tau_1$-converges to some $\mu\in\Mb(E_n)$ as soon as
$\lim_{k\rightarrow\infty}\mu_k(A)=\mu(A)$, for any $A\in\Ea_n$.
\end{defi}

We let $\theta~:~ \lim_{\Ua_n}\Ma_n\rightarrow\Mb(E_n)$ be the mapping that
associates a point $\mu=(\mu^{U_n})_{U_n\in\Ua_n}\in \lim_{\Ua_n}\Ma_n$ the set function $\theta\in \Mb(E_n)$
defined for any $A\in \Ea_n$ by
$$
\theta(\mu)(A)=\mu^{U_n}(A)\quad \mbox{\rm where $U_n\in\Ua_n$ is
such that $A\in\sigma(U_n)$}.
$$
By construction of the projective inverse spectrum and by definition of the $\tau_1$ convergence,
 it is readily checked that $\theta$ is an homeomorphism.

By Theorem~\ref{theomdp}, the random sequence $
V^N_n(U_n):=\left(V^N_n(U^1_n),\ldots,V^N_n(U^d_n)\right) $
satisfies a MDP in $\RR^d$, with speed $\alpha(N)$ and with the good
rate function
$$
I_{U_n}(v^1,\ldots,v^d):=\sup_{u\in\RR^d}{\left(\langle u,v\rangle -\frac{1}{2}~\EE\left(
 \left(\sum_{i=1}^d u^i~V_n(U^i_n)\right)^2\right)\right)}.
$$
Since we have
$$
\sum_{i=1}^d u^i~V_n(U^i_n)=V_n\left(f_u\right)\quad\mbox{\rm with}\quad
f_u:=\sum_{i=1}^d u^i~1_{U^i_n}
$$
we readily find that
$$
\frac{1}{2}~\EE\left(
 \left(\sum_{i=1}^d u^i~V_n(U^i_n)\right)^2\right)=\frac{1}{2}~\eta_{n-1}\left(
 K_{n,\eta_{n-1}}\left[f_u-K_{n,\eta_{n-1}}(f_u)\right]^2
 \right)
$$
from which we conclude that
$$
I_{U_n}(\pi_{U_n}(\mu)):=\sup_{f\in\Ba(E_n,\sigma(U_n))}{\left(\mu(f) -\frac{1}{2}~\eta_{n-1}\left(
 K_{n,\eta_{n-1}}\left[f-K_{n,\eta_{n-1}}(f)\right]^2
 \right)
\right)}.
$$

By a theorem of D. Dawson and J. Gartner, we deduce the following
\begin{prop}
The sequence of random fields
$
V^N_n
$
satisfy an MDP in $\Mb(E_n)\left(\simeq\lim_{\Ua_n}\Ma_n\right)$, with speed $\alpha(N)$ and with the good rate function
\begin{equation}\label{casdeux}
\bar I_n(\mu)=\sup_{U_n\in\Ua_n}{I_{U_n}(\pi_{U_n}(\mu))}.
\end{equation}
\end{prop}
The proof of (\ref{casun}) is now a direct consequence of the next lemma.
\begin{lem}\label{lemf}
The domain $\mbox{\rm Dom}(\bar I_n)=\left\{\mu\in\Mb(E_n)~:~\bar
I_n(\mu)<\infty\right\}$ of the mapping $\bar I_n$ is included in
$\Ma(E_n)$ and for any $\mu\in\Ma(E_n)$, the rate function $\bar
I_n(\mu)$ defined in (\ref{casdeux}) coincide with $I_n$ in
(\ref{casun}).
\end{lem}

Before getting into the proof of the lemma, it is convenient to make some remarks.

Firstly, since  the relative topology on $\Ma(E_n)$ induced by the
$\tau_1$  topology coincide with the $\tau$ topology, one concludes
that the sequence of random fields $ V^N_n $ satisfies  a MDP in
$\Ma(E_n)$ with good rate function $I_n$.

Furthermore, since the projection operators $\pi_{U_n}$ are  $\tau$-continuous,
 by the contraction principle one concludes that the random fields sequence
$ \pi_{U_n}\left(V^N_n\right) $ satisfies a MDP in
$\Ma(E_n,\sigma(U_n))$ with the good rate function
$$
I_{U_n}(\nu):=\inf{\left\{ I_n(\mu)~:~\mu\in\Ma(E_n)~\mbox{\rm
s.t.}~\pi_{U_n}(\mu)=\nu \right\}}.
$$

These constructions extend in a natural way to the sequence of random fields $(V_n^N)_{n\geq 0}$.
 Indeed, using (\ref{lapV}), we find that the random sequences
$$
\left(V^N_0(U_0),\ldots,V^N_n(U_n)\right)\quad\mbox{\rm with}\quad
(U_0,\ldots,U_n)\in(\Ua_0\times\ldots\times\Ua_n)
$$
satisfy an MDP in $\left(\RR^{d_0}\times\ldots\times\RR^{d_n}\right)$, with speed $\alpha(N)$ and with the good rate function
$$
I_{U_0,\ldots,U_n}(v_0,\ldots,v_n):=\sum_{p=0}^n\sup_{u_p\in\RR^{d_p}}{\left(
\langle u_p,v_p\rangle -\frac{1}{2}~\EE\left(
 V_p(f^{u_p}_p)^2\right)\right)}
$$
with the sequence of functions $
f_n^{u_n}=\sum_{i=1}^d u^i_n~1_{U^i_n}$.
The proof of theorem~\ref{theopr} is now easily completed.

Now, we come to the 

{\bf Proof of lemma~\ref{lemf} :} Consider a sequence of partitions
$U_{n,d}$, finer and finer when $d$ increases, such that
$\sigma\left(\bigcup_{d\ge1} U_{n,d}\right)={\mathcal E}_n$. To
prove that $\mbox{\rm Dom}(\bar I_n)\subset \Ma(E_n)$, we use the
fact that
$$
I_{U_{n,d}}(\pi_{U_{n,d}}(\mu))<\infty\Rightarrow
\pi_{U_{n,d}}(\mu)\ll \pi_{U_{n,d}}(\eta_n)
$$
and
$$
 \pi_{U_{n,d}}(\eta_n)\left(\left(\frac{d \pi_{U_{n,d}}(\mu)}{d\pi_{U_{n,d}}(\eta_n)}\right)^2\right)\leq
 I_{U_{n,d}}(\pi_{U_{n,d}}(\mu))\le \bar I_n(\mu)<\infty.
$$
See for instance (\ref{major}) in the appendix. Therefore
$\left\{\frac{d
\pi_{U_{n,d}}(\mu)}{d\pi_{U_{n,d}}(\eta_n)}\right\}_{d\ge1}$ is a
$\LL_2$-bounded martingale w.r.t. the probability measure $\eta_n$
and the filtration $(\sigma(U_{n,d}))_{d\ge1}$. By the martingale
convergence theorem, there is some $h_\mu\in\LL_2(\eta_n)$ such that
$$
\frac{d \pi_{U_{n,d}}(\mu)}{d\pi_{U_{n,d}}(\eta_n)}\to h_\mu
$$
in $\LL_2(\eta_n)$, as $d$ goes to infinity. We show now that
$h_\mu$ does not depend on the sequence $(U_{n,d})$. In fact if
$(U'_{n,d})_{d\ge1}$ is another such sequence of partitions, we
consider the  partition $V_{n,d}$ which is finer than $U_{n,d}$ and
$U'_{n,d}$ such that $V_{n, d+1}$ is finer than $V_{n,d}$. By the
above argument, we have
$$
\frac{d \pi_{U'_{n,d}}(\mu)}{d\pi_{U'_{n,d}}(\eta_n)}\to h_\mu',\
\frac{d \pi_{V_{n,d}}(\mu)}{d\pi_{V_{n,d}}(\eta_n)}\to \tilde h_\mu
$$
in $\LL_2(\eta_n)$, as $d\to\infty$. Consequently for any
$\sigma(U_{n,d})$-measurable and bounded function $f$ (with $d$
fixed),
$$
\eta_n(h_\mu f)=\eta_n\left(\frac{d
\pi_{U_{n,d}}(\mu)}{d\pi_{U_{n,d}}(\eta_n)} f\right)=
\pi_{U_{n,d}}(\mu)(f)=\pi_{V_{n,d}}(\mu)(f)=\eta_n(\tilde h_\mu f).
$$
Thus $h_\mu=\tilde h_\mu$, $\eta_n-a.s.$. By the same way
$h'_\mu=\tilde h_\mu$, $\eta_n-a.s.$. Hence $h_\mu$ does not depend
on $(U_{n,d})$.

Finally for any finite partition $U_n$ and $\sigma(U_n)$-measurable
function $f$, taking a sequence of partitions $(U_{n,d})$ containing
$U_n$, we get for $d$ large enough
$$
\mu(f)=\pi_{U_{n,d}}(\eta_n)\left(\frac{d
\pi_{U_{n,d}}(\mu)}{d\pi_{U_{n,d}}(\eta_n)} f\right)=\eta_n(fh_\mu).
$$
Consequently $\mu$ is the measure $h_\mu \eta_n$.

For the last assertion, we see that
$$\aligned
\bar I_n(\mu)&=\sup_{U_n\in\Ua_n} \sup_{f\in
\Ba(E_n,\sigma(U_n))}\left(\mu(f)- \frac 12 \EE V_n(f)^2 \right)\\
&=\sup_{f\in \bigcup_{U_n\in
\Ua_n}\Ba(E_n,\sigma(U_n))}\left(\mu(f)- \frac 12 \EE V_n(f)^2
\right)\\
&=\sup_{f\in \Ba(E_n)}\left(\mu(f)- \frac 12 \EE V_n(f)^2
\right)=I_n(\mu)
\endaligned
$$
by the fact that for any $f\in  \Ba(E_n)$, there is a sequence
$f_k\in \bigcup_{U_n\in \Ua_n}\Ba(E_n,\sigma(U_n))$ which converge
uniformly to $f$ over $E_n$, and $\EE V_n(f_k)^2\to \EE V_n(f)^2 $
by the expression of $\EE V_n(f)^2$.
 \cqfd

\subsection{Some contraction properties}

By the contraction principle, the moderate deviation principles
presented in Theorem~\ref{theopr} can be  transferred to continuous
transformations of the  local sampling random fields $V^N_n$. For
instance, we have the following proposition.

\begin{prop} \label{propcontract}
The random fields $ \sum_{p=0}^nV^N_p\Da_{p,n}$ and $ W^N_n $
satisfy the MDP in $\Ma(E_n)$ with the good rate function
\begin{eqnarray}
J_n(\nu)&=&\inf\left\{
\sum_{p=0}^nI_p(\mu_p)~:~(\mu_p)_{0\leq p\leq n}\in\prod_{p=0}^n\Ma(E_p)~\mbox{\rm s.t.}~\nu=\sum_{p=0}^n\mu_p\Da_{p,n}
\right\}\label{firstgood}\\
&=&\sup_{f\in\Ba(E_n)}{\left(
\nu(f)-\frac{1}{2}\EE\left(W_n(f)^2\right)\right)}\nonumber.
\end{eqnarray}
\end{prop}
\proof The fact that $ \sum_{p=0}^nV^N_p\Da_{p,n} $ satisfies a MDP
in $\Ma(E_n)$ with the the good rate function (\ref{firstgood}) is
an immediate consequence of theorem~\ref{theopr}. On the other hand,
using (\ref{lapV}) and (\ref{lapW})  we prove that both random
sequences
$$
W^N_n(U_n):=\left(W^N_n(U^1_n),\ldots,W^N_n(U^d_n)\right)
$$
and
$$
\sum_{p=0}^nV^N_p\Da_{p,n}(U_n):=\left(\sum_{p=0}^nV^N_p\Da_{p,n}(U^1_n),\ldots,\sum_{p=0}^nV^N_p\Da_{p,n}(U^d_n)\right)
$$
with $ U_n=(U^i_n)_{1\leq i\leq d}\in\Ua_n $, satisfies a MDP in
$\RR^d$, with speed $\alpha(N)$ and with the good rate function
$$
\Ja_{U_n}(v^1,\ldots,v^d):=
\sup_{f\in\Ba(E_n,\sigma(U_n))}{\left(\mu(f)
-\frac{1}{2}~\EE\left(W_n(f)^2\right) \right)}.
$$
We conclude that both random fields $W^N_n$ and
$\sum_{p=0}^nV^N_p\Da_{p,n}$ satisfies the same MDP in $\Ma(E_n)$
with the good rate function
$$
\Ja_{n}(\nu):=\sup_{U_n\in\Ua_n}
\sup_{f\in\Ba(E_n,\sigma(U_n))}{\left(\nu(f)
-\frac{1}{2}~\EE\left(W_n(f)^2\right) \right)}=J_n(\nu).
$$
The  the last formula comes from the uniqueness property of the rate function.
This ends the proof of the proposition.
\cqfd

\section{Moderate deviations for stochastic processes}

This section is mainly concerned with the proof of theorem~\ref{theoprpp}.
By a recent theorem~of M. A. Arcones (see for instance~theorem 3.2 in~\cite{arcones}),
this theorem is a direct consequence of the following lemma.
\begin{lem}
Under the conditions ${\bf (A1)}$ and ${\bf (A2)}$, for any $y>0$ we have
$$
\lim_{\delta\rightarrow 0}\limsup_{N\rightarrow\infty}\frac{1}{\alpha(N)}
\log{ \PP\left(\frac{1}{\sqrt{\alpha(N)}}~\left\|W_n^N\right\|_{\Fa_n(\delta)}>y\right)}=-\infty
$$
with the set of functions $\Fa_n(\delta)$ given below:
$$
\Fa_n(\delta):=\{h_n~:~h_n=(f_n-g_n)~:~(f_n,g_n)\in\Fa_n^2~:~\eta_n(h_n^2)^{1/2}\leq
\delta\}.
$$
\end{lem}

\proof
The proof of this lemma is based on several key properties of empirical processes associated with
conditionally independent sequences. These results are more or less well known, thus their are
housed in the appendix~\ref{review-empirical}.

By construction, recalling that $0\in\Fa_n$, if we choose $\delta=2$ then we have
$$
\Fa_n(\delta)=\Fa_n(2)=\left\{h=(f-g)~:~(f,g)\in \Fa_n\right\}\supset\Fa_n
$$
Thus, using elementary manipulations we prove that the condition ${\bf (A2)}$ implies that
$$
\left\|\Phi_{n+1}(\mu)-\Phi_{n+1}(\eta_{n})\right\|_{\Fa_{n+1}}
\leq c(n)~\left\|\mu-\eta_{n}\right\|_{\Sigma_n\left(\Fa_{n+1}\right)}
$$
for some separable collection $\Sigma_n\left(\Fa_{n+1}\right)$  of
measurable functions $f_n$ on $E_{n}$, s.t. $\|f_n\|\leq
1$, and such that
\begin{equation}\label{control-I}
I(\Sigma_n\left(\Fa_{n+1}\right))<c^{\prime}(n)~I(\Fa_{n+1})
\end{equation}
 for some finite constants $c(n)$ and $c^{\prime}(n)<\infty$.

 This implies that
\begin{eqnarray}
\sqrt{N}~\left\|\Phi_{n+1}(\eta^N_n)-\Phi_{n+1}(\eta_{n})\right\|_{\Fa_{n+1}}
&\leq& c(n)~\sqrt{N}\left\|\eta^N_n-\eta_{n}\right\|_{\Sigma_n\left(\Fa_{n+1}\right)}\nonumber\\
&=&c(n)~
\left\|W^N_n\right\|_{\Sigma_n\left(\Fa_{n+1}\right)}\label{Phi-Continuity}
\end{eqnarray}
On the other hand, we have
$$
W^N_{n+1}=V^N_{n+1}+\sqrt{N}~\left[\Phi_{n+1}(\eta^N_{n})-\Phi_{n+1}(\eta_{n})\right]
$$
and therefore
\begin{eqnarray*}
\left\|W^N_{n+1}\right\|_{\Fa_{n+1}}&\leq& \left\|V^N_{n+1}\right\|_{\Fa_{n+1}}+c(n)~
\left\|W^N_n\right\|_{\Sigma_n\left(\Fa_{n+1}\right)}\\
&\leq &\sum_{p=0}^{n+1} c_p(n)~
\left\|V^N_p\right\|_{\Sigma_{p,n}\left(\Fa_{n+1}\right)}
\end{eqnarray*}
with $\Sigma_{p,n}=\Sigma_{p}\circ\Sigma_{p+1,n}$, and $c_p(n)=\prod_{p\leq q<n}c(q)$.
We let
 $\pi_{\psi}[Y]$ be the Orlicz norm of an $\RR$-valued random variable $Y$
associated with the
 convex function $\psi(u)=e^{u^2}-1$, and defined by
$$
\pi_{\psi}(Y)=\inf{\{a\in (0,\infty)\;:\;\EE(\psi(|Y|/a))\leq 1\}}
$$
with the convention $\inf_{\emptyset}=\infty$. From previous calculations, we have
$$
\pi_{\psi}\left(\left\|W^N_{n+1}\right\|_{\Fa_{n+1}}\right)\leq
\sum_{p=0}^{n+1} c_p(n)~\pi_{\psi}\left(
\left\|V^N_p\right\|_{\Sigma_{p,n}\left(\Fa_{n+1}\right)}\right)
$$
Combining Lemma~\ref{lemma-Orlicz-VX} with (\ref{control-I}), we
find that
$$
\pi_{\psi}\left(\left\|W^N_{n+1}\right\|_{\Fa_{n+1}}\right)\leq
c^{\prime\prime}(n)~
I\left(\Fa_{n+1}\right)
$$
for some finite constants $c^{\prime\prime}(n)$. By (\ref{Phi-Continuity}), we also have that
$$
\sqrt{N}~\pi_{\psi}\left(
\left\|\Phi_{n+1}(\eta^N_n)-\Phi_{n+1}(\eta_{n})\right\|_{\Fa_{n+1}}
\right)\leq c^{\prime\prime\prime}(n)~
I\left(\Fa_{n+1}\right)
$$
for some finite constants $c^{\prime\prime\prime}(n)$. This shows that the random fields
$V_n^N$ satisfy the regularity condition stated in  (\ref{condition-Orlicz}).

Arguing as above, we prove that
\begin{eqnarray*}
\left\|W^N_{n}\right\|_{\Fa_{n}(\delta)}
&\leq &\sum_{p=0}^{n} \alpha_p(n)~
\left\|V^N_p\right\|_{\Fa_{p,n}(\beta_p(n)\delta)}
\end{eqnarray*}
for some separable collection $\Fa_{p,n}$  of
measurable functions $f_p$ on $E_{p}$, s.t. $\|f_p\|\leq
1$, and such that $I(\Fa_{p,n})<\infty$, and for some finite constants $ \alpha_p(n)$ and $\beta_p(n)<\infty$.

$$
 \PP\left(\left\|W_n^N\right\|_{\Fa_n(\delta)}>y\sqrt{\alpha(N)}\right)\leq \sum_{p=0}^n
  \PP\left(\left\|V^N_p\right\|_{\Ga_{p,n}(\delta)}>y_{p,n}\sqrt{\alpha(N)}\right)
$$
with $y_{p,n}={y}/{[(n+1)\alpha_p(n)]}$ and $\Ga_{p,n}(\delta):=\Fa_{p,n}(\beta_p(n)\delta)$. On the other hand, using
 lemma~\ref{key-lemma-delta}, we have
$$
\begin{array}{l}
\frac{1}{\alpha(N)}\log{  \PP\left(\left\|V^N_p\right\|_{\Ga_{p,n}(\delta)}>y_{p,n}\sqrt{\alpha(N)}\right)}\\
\\
\leq -
\frac{y_{p,n}^2}{2 a_{(\beta_p(n)\delta)}(\Fa_{p,n})^2}\left(1-\frac{\alpha(N)}{N}~y_{p,n}^2~\left(\frac{b_{(\beta_p(n)\delta)}(\Fa_{p,n})}{a_{(\beta_p(n)\delta)}(\Fa_{p,n})}\right)^2\right)\rightarrow_{N\uparrow\infty} -
\frac{y_{p,n}^2}{2 a_{(\beta_p(n)\delta)}(\Fa_{p,n})^2}
\end{array}$$
with some finite constant $ b_{\delta}(\Fa) $, and
$$
a_{\delta}(\Fa)\leq c~\int_{0}^{\delta}\sqrt{\log{\Na(\Fa,\epsilon)}}~d\epsilon \rightarrow_{\delta\downarrow 0}~0
$$
so that
$$
-\frac{y_{p,n}^2}{2 a_{(\beta_p(n)\delta)}(\Fa_{p,n})^2} \rightarrow_{\delta\downarrow 0}~-\infty
$$
This ends the proof of the lemma.
\cqfd
\section{Appendix A.}

\subsection{Proof of lemma~\ref{lem1b}}\label{seclem1b}

The first  $\LL_{m}$ almost sure estimates is a direct consequence of Kintchine's inequality,
  let us examine some direct consequences of this result.
 Combining the Lipschitz property $(\Phi_{p,n})$ of the semigroup $\Phi_{p,n}$ with the decomposition
$$
\left[\eta^N_n-\eta_n\right]=\sum_{p=0}^n~\left[ \Phi_{p,n}(\eta^N_p)- \Phi_{p,n}\left(\Phi_{p}(\eta^N_{p-1})\right) \right]
$$
 we find
that (by condition (K))
$$
\sqrt{N}~\left|\left[\eta^N_n-\eta_n\right](f_n)\right|=\sum_{p=0}^n \int~\left|V^N_p(h)\right|~T^{\Phi_{p,n}}_{\Phi_{p}(\eta^N_{p-1})}(f,dh)
$$
In the above displayed formulae, we have used the convention $\Phi_{0}(\eta^N_{-1})=\eta_0$, for $p=0$.
The proof of (\ref{eqlm}) is a direct consequence of the previous  $\LL_{m}$ almost sure estimates.
On the other hand, using decomposition
 \begin{eqnarray*}
W^N_n &=&\sqrt{N}~\sum_{p=0}^n~\left[ \Phi_{p,n}(\eta^N_p)-
\Phi_{p,n}\left(\Phi_{p}(\eta^N_{p-1})\right) \right]\end{eqnarray*}
we find that $W^N_n=I^N_n+J^N_n$ , with the pair of random measures
$(I^N_n,J^N_n)$ given by
\begin{eqnarray*}
I^N_n&:=&\sum_{p=0}^nV^N_p\Da_{p,n}^{(N)}\quad\mbox{\rm and}\quad
J^N_n:=\sqrt{N}~\sum_{p=0}^n\Ra_{p,n}\left(\eta^N_p,\Phi_{p}(\eta^N_{p-1})\right)
\end{eqnarray*}
with
$$
\Da_{p,n}^{(N)}:=\Da_{\Phi_{p}(\eta^N_{p-1})
}\Phi_{p,n}\quad\mbox{\rm and}\quad \Ra_{p,n}=\Ra^{\Phi_{p,n}}.
$$
Under our assumptions, we have the almost sure estimates
$$
\sup_{N\geq 1}\beta\left(\Da_{p,n}^{(N)}\right)\leq \beta
\left(\Da\Phi_{p,n}\right) :=\sup_{\eta\in \Pa(E_p)}\beta
\left(\Da_{\eta}\Phi_{p,n}\right).
$$
Using the generalized Minkowski integral inequality we find that
$$
N~\EE\left(\left|\Ra_{p,n}\left(\eta^N_p,\Phi_{p}(\eta^N_{p-1})\right)(f_n)\right|~\left|~
 \Aa^{(N)}_{p-1}
 \right. \right)\leq ~ \delta(R^{\Phi_{p,n}}))
$$
from which we readily conclude that
$$
 \EE\left(\left|\sqrt{N}J^N_n(f_n)\right| \right)=N
 ~\EE\left(\left|\sum_{p=0}^n\Ra_{p,n}\left(\eta^N_p,\Phi_{p}(\eta^N_{p-1})\right)(f_n)\right| \right)\leq
 ~\sum_{p=0}^n  \delta(R^{\Phi_{p,n}})$$
 The proof of (\ref{biasprop}) is now clear. This end the proof of the lemma.\cqfd

 \subsection{Proof of lemma~\ref{transferlemma}}

Using Holder inequality, for any $\delta>0$ we find that
$$
\EE\left(e^{\alpha(N)~B^N(t)}\right)\leq
\EE\left(e^{\alpha(N)~(1+\delta)A^N(t)}\right)^{\frac{1}{1+\delta}}
~\EE\left(e^{\alpha(N)\frac{1+\delta}{\delta}~\left|[A^N-B^N](t)\right|}\right)^{\frac{\delta}{1+\delta}}.
$$
Under our assumptions, we have
$$
\left| (1+\delta) A^N(t)-A^N((1+\delta)t) \right|\leq (1+\delta)~a_t(\delta)\quad\mbox{\rm with}\quad
\lim_{\delta\rightarrow 0}a_t(\delta)=0
$$
this implies that
$$
 \EE\left(e^{\alpha(N)~(1+\delta)A^N(t)}\right)^{\frac{1}{1+\delta}}\leq  \EE\left(e^{\alpha(N)~
 A^N((1+\delta)t)}\right)^{\frac{1}{1+\delta}}~e^{a_t(\delta)}
$$
and therefore
$$
\limsup_{N\rightarrow\infty}\frac{1}{\alpha(N)}\log{\EE\left(e^{\alpha(N)
(1+\delta)A^N(t)}\right)}^{\frac{1}{1+\delta}}\leq
\frac{1}{1+\delta}~\Lambda((1+\delta)t)+a_t(\delta).
$$
One conclude that
$$
\limsup_{N\rightarrow\infty}\frac{1}{\alpha(N)}\log{\EE\left(e^{\alpha(N)
B^N(t)}\right)}\leq
\frac{1}{1+\delta}~\Lambda((1+\delta)t)+a_t(\delta)\longrightarrow_{\delta\rightarrow
0} \Lambda(t).
$$
In much the same way, if we set $t^{\delta}=t/(1+\delta)$  we have
$$
\EE\left(e^{\alpha(N)A^N(t^{\delta})}\right)\leq
\EE\left(e^{\alpha(N)(1+\delta)B^N(t^{\delta})}\right)^{\frac{1}{1+\delta}}~
\EE\left(e^{\alpha(N)\frac{1+\delta}{\delta}\left|[B^N-A^N](t^{\delta})\right|}\right)^{\frac{\delta}{1+\delta}}.
$$
Under our assumptions, we have
$$
\left| (1+\delta) B^N(t/{(1+\delta)})-B^N(t) \right|\leq b_t(\delta)\quad\mbox{\rm with}\quad
\lim_{\delta\rightarrow 0}b_t(\delta)=0
$$
this implies that
$$
 \EE\left(e^{\alpha(N)~(1+\delta)B^N(t^{\delta})}\right)\leq  \EE\left(e^{\alpha(N)~B^N(t)}\right)~e^{b_t(\delta)}
$$
and therefore
$$
\liminf_{N\rightarrow\infty}\frac{1}{\alpha(N)}\log{\EE\left(e^{\alpha(N)
(1+\delta)B^N(t^{\delta})}\right)}\leq~
\liminf_{N\rightarrow\infty}\frac{1}{\alpha(N)}\log{\EE\left(e^{\alpha(N)
B^N(t)}\right)} +b_t(\delta).
$$
One conclude that
$$
(1+\delta)\Lambda(t/{(1+\delta)})\leq
\liminf_{N\rightarrow\infty}\frac{1}{\alpha(N)}\log{\EE\left(e^{\alpha(N)
B^N(t)}\right)}+b_t(\delta)
$$
and letting $\delta\downarrow 0$ we find that
$$
\Lambda(t)\leq
\liminf_{N\rightarrow\infty}\frac{1}{\alpha(N)}\log{\EE\left(e^{\alpha(N)
B^N(t)}\right)}.
$$
This ends the proof of the lemma.
\cqfd

\subsection{Proof of lemma~\ref{lem3}}\label{prooflem3}

Taking the logarithm, we find that
$$
\log{\Za^{(N)}_n(f)}=\sum_{p=0}^n\Delta_p\log{(\Za^{(N)}(f))}
$$
with the random increments
\begin{eqnarray*}
\Delta_n\log{(\Za^{(N)}(f))}
&=&\log{\EE\left(\exp{\left\{\sqrt{\alpha(N)}~V^N_n(f_n)\right\}}~|~\Aa^{(N)}_{n-1}\right)}
\end{eqnarray*}
We observe that
$$
\Delta_n\log{(\Za^N(f))}=\sum_{i=1}^N~\log{\EE\left(\exp{\left\{X^{(N,i)}_n(f_n)\right\}}~|~\Aa^{(N)}_{n-1}\right)}
$$
with the sequence of random variables
$$
X^{(N,i)}_n(f_n)= \sqrt{\frac{\alpha(N)}{N}}
\left(f_n(\xi_{n}^{N,i})-K_{n,\eta^N_{n-1}}(f_n)(\xi_{n-1}^{(N,i)})\right)
$$
such that
\begin{eqnarray*}
\EE\left(X^{(N,i)}_n(f_n)~|~\Aa^{(N)}_{n-1}\right)
&=&0\\
\EE\left(X^{(N,i)}_n(f_n)^2~|~\Aa^{(N)}_{n-1}\right)
&=&\frac{\alpha(N)}{N}~K_{n,\eta^N_{n-1}}\left[\left(f_n-K_{n,\eta^N_{n-1}}(f_n)\right)^2\right](\xi_{n-1}^{(N,i)})
\end{eqnarray*}
We recall that for every centered random variable $X$ with $|X|\leq
c$ for some $c<\infty$, we have
$$
- \frac{\sigma^2}{2}~\epsilon_1(c)~ \leq
\log{\EE(e^X)}-\frac{\sigma^2}{2}\leq~\frac{\sigma^2}{2}~
\epsilon_2(c)
$$
with the parameters $(\sigma^2,\epsilon_1(c),\epsilon_2(c))$ given
below
$$ \sigma^2=\EE(X^2)\qquad
\epsilon_1(c):=\left\{\left[1-\theta(-c)\right]+\left(\frac{\theta(-c)~c}{2}\right)^2\right\}\quad
\mbox{\rm and}\quad \epsilon_2(c):=\left[\theta(c)-1\right].
$$
In the above display,  $\theta$  is  the $\Ca^1$-increasing function
defined by $\theta(x)=\frac{2}{x^2}\left(e^x-1-x\right)$ for
$x\not=0$ and $\theta(0)=1$. We set
$$
\forall j=1,2\qquad \tau^{(N)}_{j,n}(f):=\sup_{0\leq p\leq n}
\epsilon_j\left( \sqrt{\frac{\alpha(N) }{N}  }~ \mbox{\rm osc}(f_p)
\right).
$$
Using the above estimate, for any $p\leq n$ we find that
$$
\Delta_p\log{(\Za^{(N)}(f))}-\frac{\alpha(N)}{2}~\Delta_n\langle
M^{(N)}(f)\rangle\leq \tau^{(N)}_{2,n}(f)~
\frac{\alpha(N)}{2}~\Delta_p\langle M^{(N)}(f)\rangle
$$
and
$$
-\tau^{(N)}_{1,n}(f)~\frac{\alpha(N)}{2}~\Delta_p\langle
M^{(N)}(f)\rangle\leq
\Delta_p\log{(\Za^N(f))}-\frac{\alpha(N)}{2}~\Delta_p\langle
M^{(N)}(f)\rangle.
$$
This yields that
$$-\tau^{(N)}_{1,n}(f)~\frac{\alpha(N)}{2}~\langle M^{(N)}(f)\rangle_n
\leq \log{\Za^{(N)}_n(f)}-\frac{\alpha(N)}{2}~\langle
M^{(N)}(f)\rangle_n\leq
\tau^{(N)}_{2,n}(f)~\frac{\alpha(N)}{2}~\langle M^{(N)}(f)\rangle_n.
$$
The end of the proof is now a direct consequence of the following
formula
$$
\begin{array}{l}
\exp{\left\{\sqrt{\alpha(N)}~ M^{(N)}_n(f)-\frac{\alpha(N)}{2}~\langle M^{(N)}(f)\rangle_n\right\}}\\
\\
=E^{(N)}_n(f)~\exp{\left\{
\log{\Za^{(N)}_n(f)}-\frac{\alpha(N)}{2}~\langle M^{(N)}(f)\rangle_n
\right\}}
\end{array}.
$$
This ends the proof of the lemma. \cqfd

\section{Appendix B.}\label{varform}
Given two measurable spaces $E_1$ and $E_2$, we consider  a
probability measure $\mu$ over the set $E_1$ and a Markov transition
$M(x,dy)$ from $E_1$ to $E_2$. We let $V$ be the Gaussian and
centered random field on $\LL_2(E_2,\mu M)$ defined for any
$f\in\LL_2(E_2,\mu M):=\LL_2(\mu M)$ by
\begin{equation}\label{v}
\EE\left(V(f)^2\right)=\mu\left( M([f-M(f)]^2\right).
\end{equation}
In the above display, we have used the convention
$$
x\mapsto M([f-M(f)]^2(x)=M(f^2)(x)-M(f)(x)^2\quad\mbox{\rm and}\quad
M(f)(x)=\int M(x,dy)~f(y).
$$
\begin{defi}
For any $w\in \Ma(E_2)$ we set
\begin{equation}\label{II}
I(w)=\sup_{f\in
\Ba(E)}{\left(w(f)-\frac{1}{2}~\EE\left(V(f)^2\right)\right)}.
\end{equation}
\end{defi}
\begin{lem}\label{lem1}
Assume that $E_1=E_2$ and $M(x,dy)=\mu(dy)$. In this situation, we
have
\begin{equation}\label{I}
I(w)=
\sup_{f\in \LL_2(\mu)}{\left(w(f)-\frac{1}{2}~\mu\left([f-\mu(f)]^2\right)\right)}
\end{equation}
and
\begin{equation}\label{I-iid}
I(w)=I_\mu(w):=\begin{cases}\frac 12
\left\|\frac{dw}{d\mu}\right\|_{\LL_2(\mu)}^2,\quad &\mbox{\rm
if}\quad
 w\ll\mu\quad\mbox{\rm with}\quad \frac{dw}{d\mu}\in \LL_2(\mu)\quad\mbox{and} \quad
 w(E)=0\\
 +\infty, \quad &\mbox{\rm
otherwise.}
 \end{cases}
\end{equation}

\end{lem}
\proof
 To check this claim, firstly we notice that for any constant function $f(x)=a$
we have
$$
w(f)-\mu\left([f-\mu(f)]^2\right)=a~w(E).
$$
Choosing $a=b~w(E)$ with $b\in\RR_+$, we readily find that
$$
w(E)\not=0\Rightarrow \forall b\in \RR_+\quad I(w)\geq
b~w(E)^2\Rightarrow I(w)=\infty.
$$
Whenever $w$ is not absolutely continuous w.r.t. $\mu$, we can find a measurable set $A\in\Ea$ such that $w(A)\not=0$ and $\mu(A)=0$.
 In this situation, we have
$$
\forall a\in \RR_+\quad f= a~w(A)~1_A\Rightarrow
w(f)-\mu\left([f-\mu(f)]^2\right)=a~w(A)^2\Rightarrow I(w)=\infty.
$$
On the other hand, using the fact that
$$
f=g+a\Rightarrow w(f)-\frac{1}{2}~\mu\left([f-\mu(f)]^2\right)=w(g)-\frac{1}{2}~\mu\left([g-\mu(g)]^2\right)
$$
as soon as $w(E)=0$, we can reduce the supremum in (\ref{I}) to functions $f$ with $\mu(f)=0$. This yields
$$
I(w)= \sup_{f\in
\LL_{0,2}(\mu)}{\left(\mu\left(\frac{dw}{d\mu}f\right)-\frac{1}{2}~\mu\left(f^2\right)\right)}
\quad\mbox{\rm with} \quad \LL_{0,2}(\mu)=\{h\in
\LL_2(\mu)~:~\mu(h)=0 \}.
$$
Finally, we observe that
$$
\mu\left(\frac{dw}{d\mu}f\right)-\frac{1}{2}~\mu\left(f^2\right)=\frac{1}{2}
~\mu\left(\left(f-\frac{dw}{d\mu}\right)^2\right)+\mu\left[\left(\frac{dw}{d\mu}\right)^2\right].
$$
Choosing $f=\frac{dw}{d\mu}$ we prove (\ref{I}). This ends the proof of the lemma.\cqfd

The analysis of the variational formula (\ref{II}) for more general Markov transitions $M$ is a little more involved.
Before getting into further details, we observe that
\begin{eqnarray*}
\mu\left( M([f-M(f)]^2\right)&=&(\mu M)\left( ([f-(\mu M)(f)]^2\right)-\mu\left(([M(f)-\mu M(f)]^2\right)\\
&\leq& (\mu M)\left( ([f-(\mu M)(f)]^2\right).
\end{eqnarray*}
The above inequality implies that
\begin{equation}\label{major}
I(w):=\sup_{f\in \Ba(E)}{\left(w(f)-\frac{1}{2}~\mu\left( M([f-M(f)]^2\right)\right)}\geq I_{\mu M}(w)
\end{equation}
where $I_{\mu M}(w)$ is given by (\ref{I-iid}) and therefore
$$
I(w)<\infty \Rightarrow  w\ll(\mu M)\quad\mbox{\rm with}\quad
\frac{dw}{d\mu M}\in \LL_2(\mu M)\quad\mbox{and} \quad w(E)=0.
$$
Next, we follow the
analysis developed in~\cite{liming}. Firstly, we notice that $M$ is an operator $\LL_2(\mu M)$ into $ \LL_2(\mu)$.
\begin{defi}
We let  $M^{\star}_{\mu}$ be the adjoint  operator of $M$ from
$\LL_2(\mu)$ into $ \LL_2(\mu M)$ given by
$$
\forall (f,g)\in \LL_2(\mu)\times \LL_2(\mu)\quad (\mu M)\left(f
M^{\star}_{\mu}(g)\right)=\mu (M(f)~g).
$$
\end{defi}
In fact $M_\mu^*$  can be identified as a kernel $M_\mu^*(x_2,dx_1)$
which is the conditional law of $x_1$ knowing $x_2$ under the
probability measure $\mu\otimes M (dx_1,dx_2):=\mu(dx_1)M(x_1,dx_2)$
on $E_1\times E_2$. By construction, we have
$$
\mu\left( M([f-M(f)]^2\right)=\mu(M(f^2)-M(f)^2)=(\mu M)\left(
f~(Id-M^{\star}_{\mu}M) f\right).
$$
We observe that $(M^{\star}_{\mu}M)$ is a self adjoint operator on
$\LL_2(\mu M)$ with
$$
(\mu M)=(\mu M)(M^{\star}_{\mu}M).
$$
Also notice that
$$
(\mu M)\left( f~(Id-M^{\star}_{\mu}M) f\right)=\int~\mu(dx)~M(x,dy)~\left[ f(y)-M(f)(x)\right]^2
$$
from which, we see that
\begin{equation}\label{0ex}
(\mu M)\left( f~(Id-M^{\star}_{\mu}M) f\right)=0\Leftrightarrow
f(y)=Mf(x), \ \mu(dx)~M(x,dy)-a.s.
\end{equation}
Let  $\Na$ be the subspace of those elements $h$ in $\LL_2(E_2,\mu
M)$ such that $h(y)=Mh(x), \ \mu(dx)~M(x,dy)-a.s.$. Notice that
$\Na=\{h\in L^2(E_2,\mu M); (Id-M_\mu^*M) h=0\}$. It is well known
that $\Na=\LL_2(E_2, \Ga, \mu M)$ where $\Ga$ is the
sub-$\sigma$-field generated by all $h\in \Na$ (\cite{Revuz}). In
particular $\Na \bigcap \LL_\infty(E_2,\mu M)$ is dense in $\Na$.

Consider the orthogonal supplementary subspace $\Ha_0(\mu M)$ of
$\Na$ in $\LL_2(\mu M)$.

In this notation, the rate function $I$ defined in (\ref{II}) takes
the form
$$
I(w)=\sup_{f\in \Ha_{0}(\mu M)}{\left(w(f)-\frac{1}{2}~(\mu M)\left(
f~(Id-M^{\star}_{\mu}M) f\right)\right)}.
$$

Before getting into further details, arguing as in the proof of lemma \ref{lem1}, we notice that
$$
I(w)<\infty\Rightarrow
 w\ll \mu M \quad\mbox{\rm and} \quad w(h)=0, \ \forall h\in\Na \bigcap \LL_\infty(\mu M).
$$
As $M^{\star}_{\mu}M$ is self-adjoint, definite nonnegative on
$\LL_2(\mu M)$ and its norm is $1$, we can write the spectral
decomposition of $(Id-M^{\star}_{\mu}M)$ on $\Ha_0(\mu M)$
$$
(Id-M^{\star}_{\mu}M)=\int_{[0,1]}~\lambda~dE_{\lambda}=\int_{(0,1]}~\lambda~dE_{\lambda}
$$
($E_0=0$ by the very definition of $\Ha_0(\mu M)$). The operator
$(Id-M^{\star}_{\mu}M)~:~\Ha_0(\mu M)\rightarrow \Ha_0(\mu M)$ is
injective and its inverse is given by
$$
R_{0,\mu}=(Id-M^{\star}_{\mu}M)^{-1}~:~\mbox{\rm Dom}(R_{0,\mu})\subset\Ha_0(\mu M)\mapsto
\Ha_0(\mu M)
$$
with
$$
R_{0,\mu}=\int_{(0,1]}~\frac{1}{\lambda}~dE_{\lambda}\quad \mbox{\rm
and}\quad \mbox{\rm Dom}(R_{0,\mu}):=\left\{h\in \Ha_0(\mu M) ~:~
\int_{(0,1]}~\frac{1}{\lambda^2}~d\langle
E_{\lambda}(h),h\rangle<\infty\right\}.
$$
\begin{defi}
We let $\Ha_1(\mu M)$ be the completion of the pre-Hilbert space
$\Ha_0(\mu M)$ with the inner product given by
\begin{eqnarray*}
\langle f,g\rangle_{1}&:=&\langle f,(Id-M^{\star}_{\mu}M)(g) \rangle\\
&=&\mu M\left( f~(Id-M^{\star}_{\mu}M)(g)
\right)=\int_{(0,1]}~\lambda~d \langle E_{\lambda}(f),g\rangle
\end{eqnarray*}
We define $\left(\Ha_{-1}(\mu M),\|\point\|_{-1}\right)$ as the dual space of  $\left(\Ha_{1}(\mu M),\|\point\|_{1}\right)$
w.r.t. the canonical dual relation $\Ha_0(\mu M)=\Ha_0(\mu M)^{\prime}$.
\end{defi}

By \cite{liming}, $\Ha_{-1}(E_2,\mu M)$ is the subspace of $f\in
\Ha_0(\mu M)$ such that $\|f\|_{-1}<+\infty$, and $R_{0,\mu}$ can be
regarded as an isomorphism from $\Ha_{-1}(\mu M)$ to $\Ha_{1}(\mu
M)$; furthermore for any $f\in \Ha_0(E_2,\mu M)$,
\begin{equation}
\|f\|_{-1}^2 =\int_{(0,1]}~\frac{1}{\lambda}~d\langle
E_{\lambda}(f),f\rangle=\sum_{n=0}^\infty \langle f, (M_\mu^*M)^n
f\rangle.
\end{equation}
Notice also that given $f\in \LL_2(\nu M)$ if $\sum_{n=0}^\infty
\langle f, (M_\mu^*M)^n f\rangle<+\infty$, then $f\in \Ha_0(E_2,\mu
M)$.

We further assume that $w\ll \mu M$, and $h_w=\frac{dw}{d\mu M}$
satisfies $\langle h_w, h\rangle=0$ for all $h\in \Na\bigcap
\LL_\infty(\mu M)$. In this situation, if
$$
I(w)=\sup_{f\in \Ha_{1}(\mu M)\bigcap \LL_\infty(\mu
M)}{\left(\langle h_w, f\rangle-\frac{1}{2}~\|
f\|_{1}^2\right)}<\infty
$$
then $f\to \langle h_w, f\rangle$ is a bounded linear form on
$\Ha_{1}(\mu M)\bigcap \LL_\infty(\mu M)$ w.r.t. the norm
$\|\cdot\|_1$. This yields that $h_w\in \Ha_{-1}(\mu M)$ and
\begin{eqnarray*}
\langle h_w, f\rangle-\frac{1}{2}~\| f\|_{1}^2&=&\langle R_{0,\mu}(h_w), (Id-M^{\star}_{\mu}M)(f)\rangle-\frac{1}{2}~\| f\|_{1}^2\\
&=&\langle R_{0,\mu}(h_w),f\rangle_1-\frac{1}{2}~\| f\|_{1}^2\\
&=&-\frac{1}{2}\| f-R_{0,\mu}(h_w)\|_{1}^2+\frac{1}{2}\langle R_{0,\mu}(h_w),h_w\rangle
\end{eqnarray*}
from which we conclude that
$$
I(w)=\frac{1}{2}\langle R_{0,\mu}(h_w),h_w\rangle=\frac{1}{2}\| h_w\|_{-1}^2=\frac{1}{2} \sum_{n=0}^\infty
\langle h_w, (M_\mu^*M)^n h_w\rangle.
$$
In summary we have proven

\begin{prop}\label{prop-rate} The rate function defined in (\ref{II}) is given by
$$
I(w)=\begin{cases}\frac{1}{2} \sum_{n=0}^\infty \langle h_w,
(M_\mu^*M)^n h_w\rangle, \ \ \ &\text{\rm if }\ w\ll \mu M,
h_w=\frac{dw}{d(\mu M)}\in \LL_2(\mu M)\\
+\infty, &\text{\rm otherwise.}
\end{cases}
$$

\end{prop}

\section{Appendix C.}\label{review-empirical}

In the further development of this section, $c<\infty$ stands for some finite universal constant, whose values may vary from
line to line.

Let $(\mu^i)_{i\geq 1}$
be a sequence of probability measures on a given measurable state space
$(E,\Ea)$.
During the further development of this section, we fix an integer $N
\geq 1$.
To clarify the presentation, we slightly abuse the  notation
and we denote respectively
by \index{$m(X)$}
$$
m(X)=\frac{1}{N}\sum_{i=1}^N\delta_{X^i}\quad\mbox{\rm and}\quad
\mu=\frac{1}{N}\sum_{i=1}^N\mu^i
$$ the
$N$-empirical
measure associated with a collection of independent
random variables $X=(X^i)_{i\geq 1}$, with respective distributions
$(\mu^i)_{i\geq 1}$, and the $N$-averaged measure associated with the
sequence of measures $(\mu^i)_{i\geq 1}$. We also consider the empirical random field sequences
$$
V(X)=\sqrt{N}~\left(m(X)-\mu\right)
$$
We also set
\begin{equation}\label{def-sigma-f}
\sigma(f)^2:=\EE\left(V(X)(f)^2\right)=\frac{1}{N}\sum_{i=1}^N\mu^i([f-\mu^i(f)]^2)
\end{equation}

Let $\Fa$ be a given collection  of
measurable functions $f:E\rightarrow\RR$ such that $\|f\|\leq 1$.
No generality is lost and much convenience is
gained
by supposing that the unit and the null functions
 $f=\un$ and $f=0\in\Fa$.
Furthermore, to avoid some unnecessary
technical measurability questions,
we shall also suppose that $\Fa$ is separable in the sense that it contains
a countable and dense subset. For any separable class of uniformly bounded
functions $\Ha$ s.t. $\sup_{h\in \Ha}{\|h\|}\leq H$ we set
$$
I(\Ha)=\int_{0}^{2H}\sqrt{\log{\Na(\Ha,\epsilon)}}~d\epsilon
$$

We further assume that there exists some probability measure $\overline{\mu}$ on $E$ such that
\begin{equation}\label{condition-Orlicz}
\sqrt{N}~\pi_{\psi}\left(\|\mu-\overline{\mu}\|_{\Fa}\right)\leq \tau(I(\Fa))
\end{equation}
for any class of function $\Fa$ satisfying the above properties, with finite entropy $I(\Fa)<\infty$,
and some non decreasing function $\tau$. In the above displayed formula,
 $\pi_{\psi}[Y]$ stands for the Orlicz norm of an $\RR$-valued random variable $Y$
associated with the
 convex function $\psi(u)=e^{u^2}-1$, and defined by
$$
\pi_{\psi}(Y)=\inf{\{a\in (0,\infty)\;:\;\EE(\psi(|Y|/a))\leq 1\}}
$$
with the convention $\inf_{\emptyset}=\infty$. We recall that
\begin{equation}\label{Laplace-Orlicz}
\EE\left(e^{tY}\right)\leq 2~\exp{\left(\frac{t^2}{4}~\pi_{\psi}(Y)^2\right)}
\end{equation}
for any $t\geq 0$. We prove this claim using the estimate
$$
tY=\left(\frac{t\pi_{\psi}(Y)}{\sqrt{2}}\right)~\left(\frac{\sqrt{2}~Y}{\pi_{\psi}(Y)}\right)\leq
\frac{(t\pi_{\psi}(Y))^2}{4}+\left(\frac{Y}{\pi_{\psi}(Y)}\right)^2
$$
We consider the possibly bias random field sequence
$$
\overline{V}(X)=\sqrt{N}~\left(m(X)-\overline{\mu}\right)=
V(X)+\sqrt{N}~\left(\mu-\overline{\mu}\right)
$$
The following lemma is satisfied without the regularity condition  (\ref{condition-Orlicz}).
\begin{lem}\label{lemma-Orlicz-VX}
$$
\pi_{\psi}\left(\left\|V(X)\right\|_{\Fa}\right)\leq c~I(\Fa)$$
\end{lem}
\proof
We consider a
collection of independent copies $X^{\prime}=(X^{\prime i})_{i\geq 1}$ of
the random variables $X=(X^{i})_{i\geq 1}$.
Let $\e=(\e_i)_{i\geq 1}$ constitute
a sequence that is independent and identically distributed with
$$P(\e_1=+1)=P(\e_1=-1)=1/2
$$
We also consider the empirical random field sequences
$$
V_{\epsilon}(X):=\sqrt{N}~
m_{\e}(X)
$$

We also assume that
$(\epsilon,X,X^{\prime})$ are independent. We associate with the pairs
 $(\epsilon,X)$ and $(\epsilon,X^{\prime})$  the
random measures
$
m_{\epsilon}(X)=\frac{1}{N}\sum_{i=1}^N\,\epsilon_i~\delta_{X^i}
$ and $m_{\epsilon}(X^{\prime})=\frac{1}{N}\sum_{i=1}^N\,\epsilon_i~\delta_{X^{\prime i}}$.

We notice that
\begin{eqnarray*}
\|m(X)-\mu\|^p_{\Fa}&&=\sup_{f\in\Fa}|m(X)(f)-\EE(m(X^{\prime})(f))|^p
\\
&\leq& \EE(\|m(X)-m(X^{\prime})\|^p_{\Fa}~|X)
\end{eqnarray*}
and in view of the symmetry of the random variables
$(f(X^i)-f(X^{\prime i}))_{i\geq 1}$ we have
$$
\EE(\|m(X)-m(X^{\prime})\|^p_{\Fa})=
\EE(\|m_{\epsilon}(X)-m_{\epsilon}(X^{\prime})
\|^p_{\Fa})
$$
from which we conclude that
\begin{equation}\label{VV}
E\left(\|V(X)\|_{\Fa}^p\right)\leq 2^{p}\;
E\left(\|V_{\epsilon}(X)\|_{\Fa}^p\right)
\end{equation}
By using the Chernov-Hoeffding inequality
for any $x=(x^1,\ldots,x^N)\in E^N$, the empirical process
$$
f\longrightarrow V_{\epsilon}(x)(f):=\sqrt{N}~
m_{\e}(x)(f)
$$
is sub-Gaussian for the norm
$
\|f\|_{L_2(m(x))}
=m(x)(f^2)^{1/2}
$.
Namely, for any couple of functions $f,g$ and any $\d>0$ we have
$$
\EE\left(\left[V_{\epsilon}(x)(f)-V_{\epsilon}(x)(g)\right]^2\right)=
\|f-g\|^2_{\LL_2(m(x))}
$$
and by Hoeffding's inequality
$$
P\left(\left|V_{\epsilon}(x)(f)-V_{\epsilon}(x)(g)\right|\geq \d\right)\leq 2\;e^{-\frac{1}{2}{\d^2}/{\|f-g\|^2_{\LL_2(m(x))}}}
$$
If we set $Z=\left(\frac{V_{\epsilon}(x)(f)}{\sqrt{6}\|f\|_{\LL_2(m(x))}}\right)^2$, then we find that
\begin{eqnarray*}
\EE\left(e^{Z}\right)-1&=&\int_0^{\infty} e^t~\PP\left(Z\geq t\right)~dt\\
&=&\int_0^{\infty} e^t~\PP\left(\left|V_{\epsilon}(x)(f)\right|\geq \sqrt{6t}~\|f\|_{\LL_2(m(x))}\right)~dt\\
&\leq &2~\int_0^{\infty} e^t~e^{-3t}~dt=1
\end{eqnarray*}
from which we conclude that
$$
\pi_{\psi}\left(V_{\epsilon}(x)(f)-V_{\epsilon}(x)(g)\right)\leq \sqrt{6}\|f-g\|_{\LL_2(m(x))}
$$

By the maximal inequalities for sub Gaussian processes (see
\cite{dm-ledoux}, \cite{wellner}), we find that
$$
\pi_{\psi}\left(\left\|V_{\epsilon}(x)\right\|_{\Fa}\right)\leq
c~I(\Fa)<\infty
$$
 for any $x\in E^N$. By (\ref{VV}), this clearly implies that
$$
\pi_{\psi}\left(\left\|V(X)\right\|_{\Fa}\right)\leq c~I(\Fa)$$
This ends the proof of the lemma.
\cqfd

For any $\delta>0$, we also set
$$
\Fa(\delta):=\left\{h=(f-g)~:~(f,g)\in \Fa~\mbox{\rm s.t.}~\overline{\mu}(h^2)^{1/2}\leq \delta\right\}
$$

Notice that
\begin{eqnarray*}
\Na\left(\Fa(\delta),\LL_2(\eta),\epsilon\right)&\leq&
\Na\left(\Fa(\infty),\LL_2(\eta),\epsilon\right)\leq \Na\left(\Fa,\LL_2(\eta),\epsilon/2\right)^2
\end{eqnarray*}
from which we conclude that
$$
\Na\left(\Fa(\delta),\epsilon\right)\leq \Na\left(\Fa,\epsilon/2\right)^2
$$

\begin{lem}\label{key-lemma-delta}
Under the regularity condition  (\ref{condition-Orlicz}), we have the following Laplace estimates
$$
\EE\left(e^{t\left\|V(X)\right\|_{\Fa(\delta)}
}\right)\leq 4~\exp{\left(
\frac{t^2}{2}\left[
a_{\delta}(\Fa)^2+\frac{1}{N}~\left(
tb_{\delta}(\Fa)
\right)
\right]
\right)}
$$
for any $t\geq 0$,
with the parameters
$$
a_{\delta}(\Fa)\leq c~\int_{0}^{\delta}\sqrt{\log{\Na(\Fa,\epsilon)}}~d\epsilon
$$
and
$$
b_{\delta}(\Fa)\leq c~\log{\Na(\Fa,\delta)}~\left[I(\Fa)+\tau(c~I(\Fa)\right]
$$
\end{lem}

On the other hand, for any $\delta>0$ and any $x\geq 0$, we have
$$
\log{\PP\left(\left\|V(X)\right\|_{\Fa(\delta)}\geq x\right)}\leq -\sup_{t\geq 0}{\left(
tx-\frac{t^2}{2}\left[a_{\delta}(\Fa)2+
\frac{1}{N}~
\left(tb_{\delta}(\Fa)\right)^2\right]
\right)}
$$
Explicit calculations of the Legendre-Fenchel transformation can be derived, by
choosing $t=x/\alpha_{\Fa}(\delta)^2$, we find the crude exponential concentration estimates
\begin{equation}\label{expo-delta-fin}
\frac{1}{\alpha(N)}\log{\PP\left(\frac{1}{\sqrt{\alpha(N)}}\left\|V(X)\right\|_{\Fa(\delta)}\geq x\right)}\leq -
\frac{x^2}{2 a_{\delta}(\Fa)^2}\left(1-\frac{\alpha(N)}{N}~x^2~\left(\frac{b_{\delta}(\Fa)}{a_{\delta}(\Fa)}\right)^2\right)
\end{equation}
\proof
For any probability measure $\nu$, we set
$$
d_{2,\nu}(\Fa(\delta)):=\sup_{(h_1,h_2)\in\Fa(\delta)}{\|h_1-h_2\|_{\LL_2(\nu)}}
$$
By definition, we clearly have that
\begin{equation}\label{mu2delta}
d_{2,\mu}(\Fa(\delta))\leq 2
\sup_{h\in\Fa(\delta)}{\|h\|_{\LL_2(\mu)}}\leq 2\delta
\end{equation}
Notice that for any couple of probability measures $\nu_1,\nu_2$, we have
\begin{equation}\label{nu12}
d_{2,\nu_1}(\Fa(\delta))\leq d_{2,\nu_2}(\Fa(\delta))+{\left\|\nu_1-\nu_2
\right\|^{1/2}_{\Ga(\delta)}}
\end{equation}
with
$$
\Ga(\delta)=\left\{g=(h_1-h_2)^2~:~(h_1,h_2)\in\Fa(\delta)\right\}
$$
By the maximal inequalities for sub Gaussian processes, we have the estimate
$$
\pi_{\psi}(\left\|V_{\epsilon}(x)\right\|_{\Fa(\delta)})
\leq c~\int_{0}^{d_{2,m(x)}(\Fa(\delta))}\sqrt{\log{\Na(\Fa(\delta),\epsilon)}}~d\epsilon
$$
On the other hand, using (\ref{mu2delta}) and (\ref{nu12}) we prove that the r.h.s. integral is bounded by
$$
\begin{array}{l}
 \int_{0}^{2\delta}\sqrt{\log{\Na(\Fa(\delta),\epsilon)}}~d\epsilon
\\
\\ + \int_{2\delta}^{2\delta+
 {\left\|m(x)-\overline{\mu}
\right\|^{1/2}_{\Ga(\delta)}}}\sqrt{\log{\Na(\Fa(\delta),\epsilon)}}~d\epsilon
\\
\\
\leq
 \int_{0}^{2\delta}\sqrt{\log{\Na(\Fa(\delta),\epsilon)}}~d\epsilon
 +\sqrt{\log{\Na(\Fa(\delta),2\delta)}}\times  \left\|m(x)-\overline{\mu}
\right\|^{1/2}_{\Ga(\delta)}\\
\\
\leq c~
 \int_{0}^{2\delta}\sqrt{\log{\Na(\Fa,\epsilon/2)}}~d\epsilon
 +c~\sqrt{\log{\Na(\Fa,\delta)}}\times  \left\|m(x)-\overline{\mu}
\right\|^{1/2}_{\Ga(\delta)}
\end{array}$$
We conclude that
$$
\pi_{\psi}(\left\|V_{\epsilon}(x)\right\|_{\Fa(\delta)})
\leq J_{\delta}(\Fa)+r_{\delta}(\Fa) \left\|m(x)-\overline{\mu}
\right\|^{1/2}_{\Ga(\delta)}
$$
with
$$
J_{\delta}(\Fa)\leq c~\int_{0}^{\delta}\sqrt{\log{\Na(\Fa,\epsilon)}}~d\epsilon
\quad\mbox{\rm
and}
\quad
r_{\delta}(\Fa)\leq c~\sqrt{\log{\Na(\Fa,\delta)}}
$$
Using (\ref{Laplace-Orlicz}), we have
\begin{eqnarray*}
\EE\left(e^{t\left\|V_{\epsilon}(X)\right\|_{\Fa(\delta)}
}\right)&\leq& 2~
\EE\left[e^{\frac{t^2}{2}\left(J_{\delta}(\Fa)^2+r_{\delta}(\Fa)^2 \left\|m(X)-\overline{\mu}
\right\|_{\Ga(\delta)}\right)}\right]\\
&=&2~e^{\frac{t^2}{2}J_{\delta}(\Fa)^2}~
\EE\left[e^{\frac{t^2}{2\sqrt{N}}r_{\delta}(\Fa)^2 \left\|\overline{V}(X)
\right\|_{\Ga(\delta)}
}\right]
\end{eqnarray*}
and by (\ref{VV}) we have
\begin{eqnarray*}
\EE\left(e^{t\left\|V(X)\right\|_{\Fa(\delta)}
}\right)&\leq&
\EE\left(e^{2t\left\|V_{\epsilon}(X)\right\|_{\Fa(\delta)}
}\right)\\
&\leq&2~e^{2t^2J_{\delta}(\Fa)^2}~
\EE\left[e^{\frac{2t^2}{\sqrt{N}}r_{\delta}(\Fa)^2 \left\|\overline{V}(X)
\right\|_{\Ga(\delta)}}\right]
\end{eqnarray*}
Using (\ref{Laplace-Orlicz}), we conclude that
$$
\EE\left(e^{t\left\|V(X)\right\|_{\Fa(\delta)}
}\right)\leq4\exp{\left(2t^2J_{\delta}(\Fa)^2+
\frac{(t~r_{\delta}(\Fa))^4}{N}
\pi_{\psi}\left(
\left\|\overline{V}(X)
\right\|_{\Ga(\delta)}\right)^2\right)}
$$
Our next objective is to estimate the quantity $\pi_{\psi}\left(
\left\|\overline{V}(X)
\right\|_{\Ga(\delta)}\right)$.
To this end, we let $\{h^1,\ldots,h^{n_{\epsilon/16}}\}\subset\Fa(\delta)$ be the centers of
$n_{\epsilon/16}=\Na(\Fa(\delta),\LL_2(m(x)), \epsilon/16)$ $\LL_2(x)$-balls of
radius at most
$(\epsilon/16)$ covering $\Fa(\delta)$. Using the decomposition
$$
(h_1-h_2)^2-(h^i-h^j)^2=\left[(h_1-h^i)+(h^j-h_2)\right]~
\left[(h_1-h_2)+(h^i-h^j)\right]
$$
we prove that
$$
\left|(h_1-h_2)^2-(h^i-h^j)^2\right|\leq 8
\left[|h_1-h^i|+|h^j-h_2|\right]
$$
for any $h_i,h^j\in\Fa(\delta)\left(\Rightarrow \|h_i\|\vee\|h^j\|\leq 2
\right)$. Using these estimates, we prove that
$$
\Na\left(\Ga(\delta),\LL_2(m(x)),\epsilon\right)\leq
\Na\left(\Fa(\delta),\LL_2(m(x)),\epsilon/16\right)^2\leq \Na(\Fa,\epsilon/32)^4
$$
On the other hand, we have
$$
\sup_{g\in\Ga(\delta)}{\|g\|}\leq 4 \sup_{h\in\Fa(\delta)}{\|h\|^2}\leq 16
$$
This implies that
$$
I(\Ga(\delta))=\int_{0}^{32}\sqrt{\log{\Na(\Ga(\delta),\epsilon)}}~d\epsilon\leq c~I(\Fa)
$$
and by lemma~\ref{lemma-Orlicz-VX}, we can prove that
$$
\pi_{\psi}\left(\left\|V(X)\right\|_{\Ha}\right)\leq c~~ \int_{0}^{2H}\sqrt{\log{\Na(\Ha,\epsilon)}}~d\epsilon
$$
for any class of functions $\Ha$ s.t. $\sup_{h\in \Ha}{\|h\|}\leq H$. One concludes that
$$
\pi_{\psi}\left(\left\|V(X)\right\|_{\Ga(\delta)}\right)\leq c~ I(\Fa)
$$
and therefore
$$
\pi_{\psi}\left(\left\|\overline{V}(X)\right\|_{\Ga(\delta)}\right)\leq a(\Fa):=c~ I(\Fa)+  \tau(c~I(\Fa))
$$
$$
\EE\left(e^{t\left\|V(X)\right\|_{\Fa(\delta)}
}\right)\leq4\exp{\left(2t^2J_{\delta}(\Fa)^2+
\frac{(t~r_{\delta}(\Fa))^4}{N}
a(\Fa)^2\right)}
$$
The end of the proof of the Laplace estimates is now easily completed.
This ends the proof of the theorem.
\cqfd


\begin{thebibliography}{99}

\bibitem{arcones} The large deviation principle for stochastic processes I and II. {\em Theory of
Probability and its Applications}. 47, 567-583 and 48, 19-44 (2003).
\bibitem{arcones2}
M.A. Arcones. Moderate deviations of empirical processes. Stochastic inequalities
and applications. Progr. Probab. 56, Birkhauser, Basel, 189-212 (2003).

\bibitem{deacosta}
A. de Acosta.
Moderate Deviations for Empirical Measures of Markov Chains: Lower Bounds.
{\em The Annals of Probability}, Vol. 25, No. 1, pp. 259-284  (1997).

\bibitem{deacostachen}
A. de Acosta and X. Chen. Moderate Deviations for Empirical Measures of Markov Chains: Upper Bounds. {\em
Journal of Theoretical Probability}, Vol. 11, No. 4 (1998)

\bibitem{fk}
Del Moral, P., {\em Feynman-Kac formulae. Genealogical and interacting particle systems with applications},
 Probability and its Applications, Springer Verlag, New York (2004).


\bibitem{dm-ledoux}
P. Del Moral, and M. Ledoux.
On the Convergence and the Applications of Empirical Processes for Interacting Particle Systems and Nonlinear Filtering.
{\em Journal of Theoretical Probability}, Vol. 13, No. 1, 225-257 (2000).


\bibitem{dmrio}
P. Del Moral, E. Rio, Concentration inequalities for Mean Field Particle Models.
 \href{http://hal.inria.fr/inria-00375134/en/}{\tt HAL-INRIA} publication no. 6901 (29p.), April (2009).

\bibitem{ddj} Del Moral, P., Doucet, A., Jasra, A.,
{\em Sequential Monte Carlo Samplers.}
Journal of the Royal Statistical Society, Series B, vol. 68, no. 3, pp. 411-436 (2006).

\bibitem{douc}
R. Douc, A. Guillin and J. Najim.
Moderate Deviations for Particle Filtering.
{\em The Annals of Applied Probability}, Vol. 15, No. 1B,  pp. 587-614 (2005).

\bibitem{arnaud}
Doucet A., de Freitas J.F., Gordon N.J., {\em Sequential Monte-Carlo
Methods in Practice}, Springer Verlag New York (2001).

\bibitem{gao1}
F.Q. Gao.  Moderate deviations for martingales and mixing random processes.
Stochastic Process. Appl., 61, 263?275 (1996).
\bibitem{gao2}
F.Q. Gao.  Moderate deviations and large deviations for kernel density estimators.
J. Theoret. Probab., 16, 401-418 (2003).

\bibitem{guillin1}
H. Djellout and A. Guillin.
Moderate deviations for Markov chains with atom
{\em Stochastic processes and their applications}. vol. 95, no2, pp. 203-21 (2001)

\bibitem{ledoux}
M. Ledoux, M. Talagrand. {\em Probability in Banach Spaces, Isometry and Processes}. Springer-Verlag Berlin (1991).

\bibitem{ledoux2}
M. Ledoux.  Sur les d\'eviations mod\'er\'ees de sommes de variable al\'eatoires vectorielles
ind\'ependantes de m\^eme loi. {\em Ann. Inst. Henri Poincar\'e}. 28, 267-280 (1992).

\bibitem{Revuz} D. Revuz.  {\it Markov Chains},  North-Holland , 1976.

\bibitem{talagrand}
M. Talagrand.
Sharper bounds for Gaussian and empirical processes, {\em Ann. Probab.} 22 (1994) 28-76.


\bibitem{liming2}
L. M. Wu. Large deviations, moderate deviations and LIL for
empirical processes. {\em The Annals of Probability}, vol. 22, no.
1, pp.17--27 (1994).




\bibitem{liming}
L. M. Wu, Forward-Backward martingale decomposition and compactness results for
additive functionals of stationnary ergodic Markov processes. {\em Annales de l'I.H.P,, section B, tome 35, no 2,} pp. 121--141 (1999).

\bibitem{wellner}
A. W. van der Vaart, J. A. Wellner. {\em Weak convergence of stochastic processes}. Springer Series in Statistics, Springer (1996).

\end{thebibliography}
\end{document}